\documentclass[11pt]{article}

\usepackage{amsmath,amsthm,amssymb,tikz-cd}

\newtheorem{theorem}{Theorem}[section]					
\newtheorem{lemma}[theorem]{Lemma}

\newtheorem{proposition}[theorem]{Proposition}
\theoremstyle{definition}
\newtheorem{definition}[theorem]{Definition}
\theoremstyle{definition}

\newtheorem{nothing}[theorem]{}
\newtheorem{notation}[theorem]{Notation}

\newtheorem{remark}[theorem]{Remark}

\newcommand{\K}{\mathbb{K}}
							
\newcommand{\OO}{\mathcal{O}}

\newcommand{\catname}[1]{{\normalfont\textbf{#1}}}		

\newcommand{\bimod}[2]{{}_{#1}\catname{mod}_{#2}}

\newcommand{\biset}[2]{{}_{#1}\catname{set}_{#2}}
\newcommand{\bitriv}[2]{{}_{#1}\catname{triv}_{#2}}

\newcommand{\lmod}[1]{{}_{#1}\catname{mod}}

\newcommand{\lperm}[1]{{}_{#1}\catname{perm}}
\newcommand{\lset}[1]{{}_{#1}\catname{set}}
\newcommand{\ltriv}[1]{{}_{#1}\catname{triv}}

\newcommand{\Bl}{\operatorname{Bl}}
\newcommand{\bli}{\operatorname{bli}}
\newcommand{\BP}{\mathcal{BP}}
\newcommand{\br}{\operatorname{br}}
\newcommand{\Br}{\operatorname{Br}}

\newcommand{\dsum}{\oplus}

\newcommand{\gp}[1]{\langle#1\rangle}
\newcommand{\Hom}{\operatorname{Hom}}

\newcommand{\Ind}{\operatorname{Ind}}

\newcommand{\Inj}{\operatorname{Inj}}

\newcommand{\into}{\hookrightarrow}

\newcommand{\iso}{\cong}

\newcommand{\isoto}{\overset{\sim}{\to}}

\newcommand{\lexp}[2]{\setbox0=\hbox{$#2$} \setbox1=\vbox to
                 \ht0{}\,\box1^{#1}\!#2}

\newcommand{\nor}{\trianglelefteq}

\newcommand{\onto}{\twoheadrightarrow}

\newcommand{\Res}{\operatorname{Res}}

\newcommand{\set}[1]{\left\{#1\right\}}

\newcommand{\subgp}{\leq}

\newcommand{\tensor}{\otimes}

\newcommand{\tr}{\operatorname{tr}}

	\makeatletter
	\newcommand{\tpitchfork}{%
		\vbox{
			\baselineskip\z@skip
			\lineskip-.52ex
			\lineskiplimit\maxdimen
			\m@th
			\ialign{##\crcr\hidewidth\smash{$-$}\hidewidth\crcr$\pitchfork$\crcr}
		}%
	}
	\makeatother

\title{Fixed points of extended tensor products}
\author{Robert Boltje and John Revere McHugh}

\title{Fixed points of extended tensor products\footnote{{\bf MR Subject Classification:}  
20C20, 19A22. {\bf Keywords:}  $p$-permutation module, trivial source module, $p$-permutation equivalence, isotypy, $G$-set, biset, extended tensor product}}

\author{\small Robert Boltje\\
  \small Department of Mathematics\\
  \small University of California\\
  \small Santa Cruz, CA 95064\\
  \small U.S.A.\\
  \small boltje@ucsc.edu
  \and
  \small John Revere McHugh\\
  \small Department of Mathematics\\ 
  \small University of Denver\\
  \small Denver, CO 80208\\
  \small U.S.A.\\
  \small john.r.mchugh@du.edu}
\date{November 10, 2025}

\begin{document}
\sloppy	
\maketitle

\begin{abstract}
For a $p$-permutation equivalence between two block algebras of finite groups, we introduce new square diagrams that link the $p$-permutation equivalence via the Brauer construction to local equivalences between stabilizers of corresponding Brauer pairs. These diagrams can be viewed as lifts of the square diagrams in the definition of isotypies. The proof of the commutativity requires new technical tools, namely a formula for how taking fixed points commutes with extended tensor products of finite sets with group actions and how the Brauer construction commutes with taking extended tensor products of $p$-permutation modules. These fundamental formulas, generalizing earlier results by Boltje-Danz and by Boltje-Perepelitsky, should be of independent interest.
\end{abstract}

\section{Introduction}\label{sec:1}

Throughout this article $G$, $H$, and $K$ stand for finite groups. Given subgroups $X\subgp G\times H$ and $Y\subgp H\times K$, the \textit{extended tensor product} is, roughly speaking, an operation that allows one to compose objects associated to $X$ with objects that are associated to $Y$. For example, if $R$ is a commutative ring the extended tensor product construction allows one to compose $RX$-modules with $RY$-modules. This version of the extended tensor product first appeared in \cite{Bouc_2010.1}, where it was used to provide a decomposition of tensor products of the form $\Ind_X^{G\times H}(M)\tensor_{RH}\Ind_Y^{H\times K}(N)$, for $M$ an $RX$-module and $N$ an $RY$-module. Boltje and Perepelitsky gave fundamental properties of the extended tensor product of modules in \cite[Section 6]{Boltje_2020}. This operation generalizes, for extreme choices of the subgroups $X$ and $Y$, both the tensor product of bimodules and the tensor product of representations over the ground ring $R$. For example, when $X=G\times H$ and $Y=H\times K$ the extended tensor product of an $RX$-module $M$ with an $RY$-module $N$ is nothing other than the tensor product $M\tensor_{RH}N$.

Here we introduce the extended tensor product of finite $X$-sets and $Y$-sets (see Definition \ref{defn:exttensprod} below). This version of the construction can be viewed as a functor 
\begin{equation*}
	-\underset{X,Y}{\hat{\tensor}}-:\lset{X}\times\lset{Y}\to\lset{X\ast Y}
\end{equation*}
where $X\ast Y\subgp G\times K$ denotes the composition of $X$ and $Y$ (cf. \cite[2.3.19]{Bouc_2010} or \ref{comp} below). If $R$ is a commutative ring then, via the linearization functor, the extended tensor product of finite $X$-sets and $Y$-sets commutes with the extended tensor product of $RX$-modules and $RY$-modules.

Now, according to Burnside's theory of marks, any finite left $G$-set $U$ is determined up to isomorphism by the $N_G(H)$-subsets $U^H$ of $H$-fixed points as $H$ runs over all subgroups of $G$. Therefore, if $U$ is a finite left $X$-set and $V$ is a finite left $Y$-set, it is desirable to have information about the fixed point subsets $(U\hat{\tensor}_{X,Y}V)^Z$, where $Z\subgp X\ast Y$. The main theorem of this article yields a decomposition of such subsets, provided a condition on the point-stabilizers is satisfied:

\begin{theorem}\label{thm:1}
	Let $X\subgp G\times H$, $Y\subgp H\times K$, and $Z\subgp X\ast Y$. Let $U\in\lset{X}$, $V\in\lset{Y}$ and suppose that $k(X_u,Y_v)=1$ for all $u\in U$ and $v\in V$. Let $[T\backslash\Omega]$ denote a set of representatives for the orbits of $T$ on $\Omega$. Then there is a natural isomorphism of $N_{X\ast Y}(Z)$-sets
	\begin{equation*}
		\eta_{U,V}:\bigsqcup_{\omega\in[T\backslash\Omega]}\Ind_{N_{X\ast Y}^\omega(Z)}^{N_{X\ast Y}(Z)}
		\Bigl(\Res_{N_{X\ast Y}^\omega(Z)}^{\tilde{X}(\omega)\ast \tilde{Y}(\omega)}
		\bigl(U^{X(\omega)}\underset{\tilde{X}(\omega),\tilde{Y}(\omega)}{\hat{\tensor}}V^{Y(\omega)}\bigr)\Bigr)
		\isoto (U\underset{X,Y}{\hat{\tensor}}V)^Z
	\end{equation*}
	defined by
	\begin{equation*}
		\eta_{U,V}\bigl((g,k)\tensor(u\tensor v)\bigr)=(g,k)\cdot(u\tensor v)
	\end{equation*}
	for all $\omega\in[T\backslash\Omega]$, $(g,k)\in N_{X\ast Y}(Z)$, $u\in U^{X(\omega)}$, and $v\in V^{Y(\omega)}$.
\end{theorem}

We note that Theorem \ref{thm:1} is a generalization of \cite[Theorem 2.3]{Boltje_2012}, where the question of decomposing fixed point subsets of the (usual) tensor product of two left-free bisets is considered.

The statement of Theorem \ref{thm:1} above does not make sense without definitions of $T$, $\Omega$, $X(\omega)$, etc. We provide these in Section \ref{sec:setupthm1} before the proof of Theorem \ref{thm:1} in Section \ref{sec:pfthm1}. For the purposes of this introduction, we only wish to observe that Theorem \ref{thm:1} provides a decomposition of the $N_{X\ast Y}(Z)$-set $(U\hat{\tensor}_{X,Y}V)^Z$ into a disjoint union of $N_{X\ast Y}(Z)$-sets which are each induced from extended tensor products of fixed point subsets. In other words, the theorem details how the extended tensor product commutes with taking fixed points. We remark that even if one sought only to determine how the (usual) tensor product of bisets commutes with taking fixed points, one would ultimately be led towards the extended tensor product operation because if $U$ is a $(G,H)$-biset, $V$ is an $(H,K)$-biset, $X\subgp G\times H$, and $Y\subgp H\times K$, then $U^X$ is a left $N_{G\times H}(X)$-set and $V^Y$ is a left $N_{H\times K}(Y)$-set, and it is not obvious how to ``compose'' $U^X$ and $V^Y$, since the subgroups $N_{G\times H}(X)$ and $N_{H\times K}(Y)$ need not be sub-direct products. 

In Section \ref{sec:apps} we present two applications of Theorem \ref{thm:1}. In fact, these applications are the motivation behind Theorem~\ref{thm:1}. The first application, Theorem \ref{thm:2}, essentially results from ``linearizing'' Theorem \ref{thm:1} over a field $F$ of positive characteristic $p$, and provides the corresponding statement for direct summands of permutation modules defined over $F$. These modules are known in the literature as $p$\textit{-permutation} modules, and also as \textit{trivial source} modules. In this situation, the \textit{Brauer construction} --- an operation defined for each $p$-subgroup of the group considered ---- replaces the operation of ``taking fixed points.'' The relevant definitions and facts are recalled in Section \ref{sec:trivsourcemodsandbrauerpairs}. We note that our Theorem \ref{thm:2} generalizes both \cite[Corollary 7.4(b)]{Boltje_2020} and \cite[Theorem 3.3]{Boltje_2012}.

Finally, in Theorem \ref{thm:3} we construct new commutative diagrams associated to a $p$-permutation equivalence between blocks of finite group algebras. Let $G$ and $H$ be finite groups, let $F$ be a field of positive characteristic $p$, which is large enough for $G$ and $H$, let $A$ be a block algebra of $FG$, and let $B$ be a block algebra of $FH$. Recall that a virtual Morita equivalence $\gamma$ between $A$ and $B$ is a $p$\textit{-permutation equivalence} if the indecomposable $F[G\times H]$-modules that support $\gamma$ are all $p$-permutation modules that have twisted diagonal vertices (cf. \cite{Boltje_2020} or \ref{ppermequivs} below). If $\gamma$ is a $p$-permutation equivalence between $A$ and $B$, Boltje and Perepelitsky showed \cite[Theorems 10.11, 11.2]{Boltje_2020} that $\gamma$ selects a maximal $A$-Brauer pair $(D,e_D)$, a maximal $B$-Brauer pair $(E,f_E)$, and an isomorphism $\phi:E\isoto D$ that induces an isomorphism between the fusion systems of the blocks $A$ and $B$. We show that if $(Q,f_Q)$ is a $B$-Brauer pair with corresponding $A$-Brauer pair $(P,e_P)$, then $\gamma$ induces commutative diagrams between the Grothendieck groups of $p$-permutation modules belonging to $A$, $B$, $FN_G(P,e_P)e_P$, and $FN_H(Q,f_Q)f_Q$. In Remark \ref{rmk:isotypydiagrams}, we show how these diagrams can be considered as liftings of the commutative diagrams that appear in Brou\'{e}'s concept of an \textit{isotypy} between blocks (cf. \cite{Broue_1990}) to the level of trivial source modules. One novel feature is that the diagrams of Theorem \ref{thm:3} are defined for every subgroup $Q$ of the defect group $E$ of $B$, rather than for every \textit{cyclic} subgroup of $E$, as in the case of an isotypy.

\section{Preliminaries and the Extended Tensor Product}

\begin{nothing}\textit{Group theoretic notation.}
	We use standard group-theoretic notation throughout this note. In particular, if $g,h\in G$ then ${}^gh=ghg^{-1}$ and $h^g=g^{-1}hg$. If $S\subseteq G$ then ${}^gS=\set{{}^gh|h\in S}$ and $c_g:S\to {}^gS$ denotes the map defined by $c_g(h)={}^gh$ for all $h\in S$. If $H\subgp G$ then $N_G(H)$ and $C_G(H)$ denote, respectively, the normalizer and centralizer of $H$ in $G$.
\end{nothing}

\begin{nothing}$G$\textit{-sets and }$(G,H)$\textit{-bisets.} (cf. \cite[Chapter 2]{Bouc_2010})
	The categories of finite left $G$-sets and finite $(G,H)$-bisets are denoted by $\lset{G}$ and $\biset{G}{H}$, respectively. Recall that any $(G,H)$-biset $U$ may be regarded as a left $G\times H$-set (and vice versa) by defining $(g,h)u=guh^{-1}$ for all $(g,h)\in G\times H$ and $u\in U$. ``Switching sides'' in this manner defines an isomorphism of categories $\biset{G}{H}\iso\lset{G\times H}$. 
	
	If $U$ is a left (respectively, right) $G$-set then $G\backslash U$ (resp., $U/G$) denotes the set of $G$-orbits of elements of $U$. If $u\in U$ then $G_u$ denotes the stabilizer in $G$ of $u$ and $[u]$ denotes the orbit of $u$. If $H\subgp G$ then $U^H$ denotes the subset of $H$-fixed points in $U$. Recall that $U^H$ inherits an action of the normalizer $N_G(H)$ of $H$ in $G$.
\end{nothing}

\begin{nothing}\label{tensprodbiset}\textit{The tensor product of bisets.} (cf. \cite[Definition 2.3.11]{Bouc_2010})
	If $U\in\biset{G}{H}$ and $V\in\biset{H}{K}$ then the \textit{tensor product} (originally known as the \textit{composition}) of $U$ and $V$ over $H$ is the $(G,K)$-biset $U\tensor_H V$ defined as follows: the elements of $U\tensor_H V$ are the $H$-orbits of $U\times V$ under the left action $h\cdot(u,v)=(uh^{-1},hv)$. The $H$-orbit of $(u,v)$ is denoted by $u\tensor v$. Thus, every element of $U\tensor_H V$ is a ``simple tensor.'' The $(G,K)$-biset structure comes from defining $g\cdot (u\tensor v)\cdot k=gu\tensor vk$ for all $g\in G$, $u\in U$, $v\in V$, and $k\in K$. 
\end{nothing}

\begin{nothing}\label{induction}\textit{Induction.} (cf. \cite[2.3.9]{Bouc_2010})
	Let $G$ be a finite group, $H\subgp G$. Via left and right multiplication we may consider $G$ as a $(G,H)$-biset. If $U\in\lset{H}$ then $\Ind_H^G(U)$ is by definition the left $G$-set $G\tensor_H U$. For any $V\in\lset{G}$ one has a natural isomorphism
	\begin{equation*}
		\Hom_G(\Ind_H^G(U),V)\iso\Hom_H(U,\Res_H^G(V)),
	\end{equation*}
	where $\Res_H^G(V)$ denotes the $H$-set obtained by restricting the action of $G$ on $V$ to $H$. If $\varphi:U\to\Res_H^G(V)$ is an $H$-equivariant map then the corresponding $G$-equivariant map $\hat{\varphi}:\Ind_H^G(U)\to V$ is defined by $\hat{\varphi}(g\tensor u)=g\varphi(u)$ for all $g\in G$, $u\in U$. 
\end{nothing}

\begin{lemma}\label{lem:injmapfrominduct}
	Let $G$ be a finite group, $H\subgp G$, $U\in\lset{H}$ and $V\in\lset{G}$. If $\varphi:U\to\Res_H^G(V)$ is an $H$-equivariant map satisfying
	\begin{itemize}
		\item[(1)] $\varphi$ is injective; and
		\item[(2)] Whenever $g\in G$ and $u,u'\in U$ are such that $g\varphi(u)=\varphi(u')$ then $g\in H$,
	\end{itemize}
	then the corresponding $G$-equivariant map $\hat{\varphi}:\Ind_H^G(U)\to V$ is injective.
\end{lemma}

\begin{proof}
	Let $\varphi:U\to\Res_H^G(V)$ be such a map and suppose that $g,g'\in G$ and $u,u'\in U$ are such that $\hat{\varphi}(g\tensor u)=\hat{\varphi}(g'\tensor u')$. Then $g\varphi(u)=g'\varphi(u')$, so $g^{-1}g'\in H$ by condition (2). Since $\varphi$ is $H$-equivariant we have $\varphi(u)=g^{-1}g'\varphi(u')=\varphi(g^{-1}g'u')$, and since $\varphi$ is assumed injective it follows that $u=g^{-1}g'u'$. Finally, within $\Ind_H^G(U)$ we have
	\begin{equation*}
		g\tensor u=g(g^{-1}g')\tensor(g^{-1}g')^{-1}u=g'\tensor u',
	\end{equation*}
	so $\hat{\varphi}$ is injective.
\end{proof}

\begin{nothing}\label{subgpdirprod}\textit{Subgroups of direct products.} (cf. \cite[2.3.21]{Bouc_2010})
	Whenever a direct product $G\times H$ of finite groups is formed, write $p_1:G\times H\onto G$ and $p_2:G\times H\onto H$ for the canonical projections. If $X\subgp G\times H$ set
	\begin{equation*}
		k_1(X)=\set{g\in G|(g,1)\in X}\qquad\text{and}\qquad k_2(X)=\set{h\in H|(1,h)\in X}.
	\end{equation*}
	Then $k_1(X)\nor p_1(X)$, $k_2(X)\nor p_2(X)$, and $k_1(X)\times k_2(X)\nor X$. Moreover, the projection maps $p_i$ induce group isomorphisms 
	\begin{equation*}
		\frac{p_1(X)}{k_1(X)}\overset{\sim}{\longleftarrow}\frac{X}{k_1(X)\times k_2(X)}\overset{\sim}{\longrightarrow}\frac{p_2(X)}{k_2(X)}.
	\end{equation*}
\end{nothing}

\begin{nothing}\label{comp}\textit{Composition of subgroups.} (cf. \cite[2.3.19]{Bouc_2010})
	If $X\subgp G\times H$ and $Y\subgp H\times K$ then the \textit{composition} of $X$ and $Y$ is the subgroup $X\ast Y$ of $G\times K$ consisting of pairs $(g,k)$ for which there exists an element $h\in H$ satisfying $(g,h)\in X$ and $(h,k)\in Y$. By \cite[Lemma 2.3.22]{Bouc_2010} we always have 
	\begin{equation*}
		k_1(X)\subgp k_1(X\ast Y)\subgp p_1(X\ast Y)\subgp p_1(X)
	\end{equation*}
	and
	\begin{equation*}
		k_2(Y)\subgp k_2(X\ast Y)\subgp p_2(X\ast Y)\subgp p_2(Y).
	\end{equation*}
	In particular, $k_1(X)\times k_2(Y)\nor X\ast Y$. Note also that if $X'\subgp X$ and $Y'\subgp Y$, then $X'\ast Y'\subgp X\ast Y$.
\end{nothing}

\begin{nothing}\label{twisteddiagsubg}\textit{Twisted diagonal subgroups.}
	Subgroups of $G\times H$ of the form $\Delta(P,\phi,Q)=\set{(\phi(y),y)|y\in Q}$, where $P\subgp G$, $Q\subgp H$ and $\phi:Q\isoto P$ is an isomorphism, are called \textit{twisted diagonal subgroups}. We may also write $\Delta_\phi(Q)$ in place of $\Delta(P,\phi,Q)$. If $P\subgp G$, $Q\subgp H$ and $\psi:P\isoto Q$ is an isomorphism then $\Delta^\psi(P)=\set{(x,\psi(x))|x\in P}$. Therefore $\Delta(P,\phi,Q)=\Delta_\phi(Q)=\Delta^{\phi^{-1}}(P)$. Note that a subgroup $X$ of $G\times H$ is twisted diagonal if and only if $k_1(X)$ and $k_2(X)$ are both trivial. The composition of two twisted diagonal subgroups, when defined, is twisted diagonal. We have the following conjugation formula for twisted diagonal subgroups of $G\times H$:
	\begin{equation*}
		{}^{(g,h)}\Delta(P,\phi,Q)=\Delta({}^gP,c_g\phi c_h^{-1},{}^hQ).
	\end{equation*}
\end{nothing}

\begin{nothing}\label{linkelts}\textit{Linking elements.}
	Let $X\subgp G\times H$ and $Y\subgp H\times K$. For each $(g,k)\in G\times K$ we define
	\begin{equation*}
		\ell_{X,Y}(g,k)=\set{h\in H|(g,h)\in X\text{ and }(h,k)\in Y}.
	\end{equation*}
	Elements of $\ell_{X,Y}(g,k)$ are said to \textit{link} the elements $g$ and $k$ or are simply called \textit{linking elements}. If the subgroups $X$ and $Y$ are contextually clear we may write $\ell(g,k)$ in place of $\ell_{X,Y}(g,k)$. Note that the composition $X\ast Y$ consists precisely of those ordered pairs $(g,k)\in G\times K$ for which $\ell(g,k)$ is nonempty.
\end{nothing}

\begin{nothing}\label{subgpskppz}\textit{The subgroups }$k(X,Y)$\textit{, }$p(X,Y)$\textit{, and }$p_Z(X,Y)$.
	If $X\subgp G\times H$ and $Y\subgp H\times K$ then set
	\begin{equation*}
		k(X,Y)=k_2(X)\cap k_1(Y)\subgp H\qquad\text{and}\qquad p(X,Y)=p_2(X)\cap p_1(Y)\subgp H,
	\end{equation*}
	where $k_2(X)$, $k_1(Y)$, etc. are as in \ref{subgpdirprod} above. Note that $k(X,Y)\nor p(X,Y)$. It is straightforward to check that if $(g,k)\in X\ast Y$ and $h\in\ell_{X,Y}(g,k)$ then $h\in p(X,Y)$ and $\ell_{X,Y}(g,k)=h\cdot k(X,Y)$. In fact, $\ell_{X,Y}$ defines a surjective group homomorphism $\ell_{X,Y}:X\ast Y\onto p(X,Y)/k(X,Y)$ with kernel $k_1(X)\times k_2(Y)$ \cite[Lemma 2.3.22]{Bouc_2010}. For this reason we set
	\begin{equation*}
		\ell(X,Y)=p(X,Y)/k(X,Y).
	\end{equation*} 
	
	
	More generally, if $Z\subgp X\ast Y$ let $p_Z(X,Y)$ denote the subgroup of elements $h\in H$ for which there exists a pair $(g,k)\in Z$ such that $h\in\ell_{X,Y}(g,k)$. Then $k(X,Y)\nor p_Z(X,Y)\subgp p(X,Y)$. We set $\ell_Z(X,Y)=p_Z(X,Y)/k(X,Y)$, and note that $\ell_Z(X,Y)$ is equal to the image of $Z$ under the homomorphism $\ell_{X,Y}$.
\end{nothing}

\begin{definition}\label{defn:exttensprod}
	Let $X\subgp G\times H$, $Y\subgp H\times K$, $U\in\lset{X}$ and $V\in\lset{Y}$. The \textit{extended tensor product} of $U$ and $V$ is the left $X\ast Y$-set $U\hat{\tensor}_{X,Y}V$ defined as follows: the elements of $U\hat{\tensor}_{X,Y}V$ are the $k(X,Y)$-orbits of $U\times V$ under the action $h\cdot(u,v)=((1,h)u,(h,1)v)$. The $k(X,Y)$-orbit of $(u,v)$ is denoted by $u\tensor v$. The action of $X\ast Y$ on $U\hat{\tensor}_{X,Y}V$ is defined by
	\begin{equation*}
		(g,k)\cdot(u\tensor v)=(g,h)u\tensor(h,k)v
	\end{equation*}
	where $(g,k)\in X\ast Y$, $u\in U$, $v\in V$ and $h\in H$ is any element linking $g$ and $k$. When $X$ and $Y$ are contextually clear we may write $U\hat{\tensor} V$ in place of $U\hat{\tensor}_{X,Y}V$. 
\end{definition}

\begin{nothing}\label{rmksondefnext}\textit{Remarks on Definition \ref{defn:exttensprod}.}
	Continuing with the notation above, we remark that:\\
	\textbf{(a)} Every element of $U\hat{\tensor}_{X,Y}V$ is a ``simple tensor'' of the form $u\tensor v$ where $u\in U$ and $v\in V$ (just as for the tensor product of bisets).\\
	\textbf{(b)} The action of $X\ast Y$ on $U\hat{\tensor}_{X,Y}V$ is well-defined: indeed, if $(g,k)\in X\ast Y$, $u\in U$, and $v\in V$ then $(g,k)\cdot(u\tensor v)$ does not depend on the choice of linking element, since if $h,h'\in\ell_{X,Y}(g,k)$ then $h'h^{-1}\in k(X,Y)$ and
	\begin{align*}
		(g,h)u\tensor(h,k)v	&=(1,h'h^{-1})(g,h)u\tensor(h'h^{-1},1)(h,k)v\\
								&=(g,h')u\tensor(h',k)v.
	\end{align*}
	Also, if $u'\in U$ and $v'\in V$ are such that $u'\tensor v'=u\tensor v$ then there exists an element $h'\in k(X,Y)$ such that $(1,h')u=u'$ and $(h',1)v=v'$, and therefore if $h\in\ell_{X,Y}(g,k)$ then we have
	\begin{align*}
		(g,k)\cdot(u'\tensor v')	&=(g,h)u'\tensor(h,k)v'=(g,hh')u\tensor(hh',k)v\\
									&=(g,k)\cdot(u\tensor v).
	\end{align*}

	\noindent\textbf{(c)} It follows immediately from the definitions that
	\begin{equation*}
		\Res_{k_1(X)\times k_2(Y)}^{X\ast Y}(U\underset{X,Y}{\hat{\tensor}} V) \iso \bigl(\Res_{k_1(X)\times k(X,Y)}^X U\bigr)\underset{k(X,Y)}{\tensor}\bigl(\Res_{k(X,Y)\times k_2(Y)}^YV\bigr).
	\end{equation*}
	Note that in the case where $X=G\times H$ and $Y=H\times K$ one has $X\ast Y=G\times K$ and the isomorphism above simply states that the extended tensor product of a left $G\times H$-set $U$ and a left $H\times K$-set $V$ is isomorphic to the tensor product of $U$ and $V$ over $H$. In other words, the extended tensor product operation generalizes the tensor product of bisets.\\
	\textbf{(d)} The extended tensor product is characterized by the universal property of Proposition \ref{prop:univpropexttens} below, from which it is straightforward to check that the extended tensor product defines a functor
	\begin{equation*}
		-\underset{X,Y}{\hat{\tensor}}-:\lset{X}\times\lset{Y}\to\lset{X\ast Y}\,,
	\end{equation*}
which is {\em biadditive}, i.e., it commutes with coproducts in each arguments.\\
	\textbf{(e)} The notation we have chosen for the extended tensor product is different than that employed in \cite{Boltje_2020}, where the operation is written $-\underset{H}{\overset{X,Y}{\tensor}}-$. Some reasons for the change of notation are: (1) since the extended tensor product is defined relative to the subgroups $X$ and $Y$, we prefer to write $X,Y$ underneath the tensor symbol rather than above it, (2) it is almost always clear from the context what the group $H$ is, so the symbol for $H$ can be safely omitted, and (3) when it is convenient, and when $X$ and $Y$ are contextually clear, we would like to drop the symbols $X$ and $Y$ from the notation, and including a hat symbol above the tensor will remind us that the extended tensor product is meant.
\end{nothing}

\begin{proposition}\label{prop:univpropexttens}
	Let $X\subgp G\times H$, $Y\subgp H\times K$, $U\in\lset{X}$ and $V\in\lset{Y}$. If $W\in\lset{X\ast Y}$ and $\varphi:U\times V\to W$ is a map satisfying
	\begin{equation*}
		(g,k)\varphi(u,v)=\varphi\bigl((g,h)u,(h,k)v\bigr)
	\end{equation*}
	for all $(g,k)\in X\ast Y$, all $h\in\ell_{X,Y}(g,k)$, and all $(u,v)\in U\times V$ then there exists a unique $X\ast Y$-equivariant map $\overline{\varphi}:U\hat{\tensor}_{X,Y}V\to W$ making the diagram
	\[\begin{tikzcd}
		{U\times V} & W \\
		{U\hat{\tensor}_{X,Y}V}
		\arrow["\varphi", from=1-1, to=1-2]
		\arrow["{-\tensor-}"', two heads, from=1-1, to=2-1]
		\arrow["{\overline{\varphi}}"', dashed, from=2-1, to=1-2]
	\end{tikzcd}\]
	commute. For $(u,v)\in U\times V$ one has $\overline{\varphi}(u\otimes v)=\varphi(u,v)$.
\end{proposition}

\begin{proof}
	Let $(u,v)\in U\times V$ and let $h\in k(X,Y)$. Then $h\in\ell_{X,Y}(1,1)$, so by assumption we have
	\begin{equation*}
		\varphi\bigl((1,h)u,(h,1)v\bigr)=(1,1)\varphi(u,v)=\varphi(u,v).
	\end{equation*}
	It follows that $\varphi$ induces a well-defined map $\overline{\varphi}:U\hat{\tensor} V\to W$, $u\tensor v\mapsto\varphi(u,v)$. The map $\overline{\varphi}$ is $X\ast Y$-equivariant, since if $(g,k)\in X\ast Y$ and $h\in\ell_{X,Y}(g,k)$ then
	\begin{align*}
		(g,k)\overline{\varphi}(u\tensor v)	&=(g,k)\varphi(u,v)=\varphi\bigl((g,h)u,(h,k)v\bigr)=\overline{\varphi}\bigl((g,h)u\tensor(h,k)v\bigr)\\
		&=\overline{\varphi}\bigl((g,k)\cdot(u\tensor v)\bigr)
	\end{align*}
	for any $u\in U$ and $v\in V$. It is clear that $\overline{\varphi}$ makes the diagram in question commute and that it is the unique $X\ast Y$-equivariant map with this property.
\end{proof}

\begin{lemma}\label{lem:stabincomp}
	Let $X\subgp G\times H$, $Y\subgp H\times K$, $U\in\lset{X}$ and $V\in\lset{Y}$. If $u\in U$ and $v\in V$ then
	\begin{equation*}
		(X\ast Y)_{u\tensor v}=X_u\ast Y_v.
	\end{equation*}
\end{lemma}

\begin{proof}
	Let $(g,k)\in (X\ast Y)_{u\tensor v}$, so that $(g,k)\cdot(u\tensor v)=u\tensor v$. By definition of the $X\ast Y$-action on $U\hat{\tensor}_{X,Y}V$ we have $(g,h)u\tensor(h,k)v=u\tensor v$ for any $h\in\ell_{X,Y}(g,k)$. Keeping such an $h$ fixed, the pairs $\bigl((g,h)u,(h,k)v\bigr)$ and $(u,v)$ must lie in the same $k(X,Y)$-orbit. In other words, there exists an element $h'\in k(X,Y)$ such that $(g,h'h)u=u$ and $(h'h,k)v=v$. Then $h'h\in\ell_{X_u,Y_v}(g,k)$, hence $(g,k)\in X_u\ast Y_v$ and since $(g,k)$ was arbitrary we have $(X\ast Y)_{u\tensor v}\subgp X_u\ast Y_v$.
	
	Suppose conversely that $(g,k)\in X_u\ast Y_v$. Since $X_u\subgp X$ and $Y_v\subgp Y$ we have $X_u\ast Y_v\subgp X\ast Y$, so $(g,k)\in X\ast Y$. By definition of the composition of subgroups there exists an element $h\in H$ such that $(g,h)\in X_u$ and $(h,k)\in Y_v$. We have $(g,k)\cdot(u\tensor v)=(g,h)u\tensor(h,k)v=u\tensor v$, hence $(g,k)\in(X\ast Y)_{u\tensor v}$ and since $(g,k)$ was arbitrary we conclude that $X_u\ast Y_v\subgp(X\ast Y)_{u\tensor v}$. The proof is complete.
\end{proof}

\section{Setup for Theorem \ref{thm:1}}\label{sec:setupthm1}

In this section we set up all of the definitions and notations needed to make the statement of Theorem \ref{thm:1} precise. We also establish key lemmas that will be used in (and help to motivate the ideas behind) the proof of Theorem \ref{thm:1}. 

We continue to let $G$, $H$, and $K$ stand for finite groups. Throughout this section let $X\subgp G\times H$, $Y\subgp H\times K$, and $Z\subgp X\ast Y$.

Given a finite left $X$-set $U$ and a finite left $Y$-set $V$, our goal is to provide a decomposition of the $N_{X\ast Y}(Z)$-set $(U\hat{\tensor}_{X,Y}V)^Z$. We will be able to accomplish this so long as we make the additional assumption that $k(X_u,Y_v)=1$ for all $u\in U$ and $v\in V$. The next lemma indicates why this assumption is needed and motivates the definitions that follow its proof.

\begin{lemma}\label{foretasurj}
	Let $U\in\lset{X}$ and $V\in\lset{Y}$. Suppose that $u\in U$ and $v\in V$ are such that $u\tensor v\in(U\hat{\tensor}_{X,Y}V)^Z$ and $k(X_u,Y_v)=1$. Then for any $(g,k)\in Z$ there exists a unique element $\omega(g,k)\in H$ satisfying
	\begin{equation*}
		(g,\omega(g,k))\in X_u\qquad\text{and}\qquad(\omega(g,k),k)\in Y_v.
	\end{equation*}
	Moreover, $\omega$ defines a homomorphism $Z\to p_Z(X,Y)$. 
\end{lemma}

\begin{proof}
	Let $(g,k)\in Z$. Then $(g,k)\in (X\ast Y)_{u\tensor v}$. Since $(X\ast Y)_{u\tensor v}=X_u\ast Y_v$ by Lemma \ref{lem:stabincomp}, there exists an element $h\in H$ such that $(g,h)\in X_u$ and $(h,k)\in Y_v$, i.e., $h\in\ell_{X_u,Y_v}(g,k)$. Now $k(X_u,Y_v)=1$ by assumption, so by the remarks in \ref{subgpskppz} above we have $\ell_{X_u,Y_v}(g,k)=h\cdot k(X_u,Y_v)=\set{h}$. Thus if we set $\omega(g,k)=h$ then $\omega(g,k)$ is the unique element of $H$ satisfying $(g,\omega(g,k))\in X_u$ and $(\omega(g,k),k)\in Y_v$. Note that $\omega(g,k)$ belongs to both $\ell_{X,Y}(g,k)$ and $p_Z(X,Y)$.
	
	Since $k(X_u,Y_v)=1$ the group $\ell(X_u,Y_v)$ defined in \ref{subgpskppz} is canonically isomorphic to $p(X_u,Y_v)$. It is clear that $p(X_u,Y_v)\subgp p(X,Y)$ and that $\omega$ is the composite of the homomorphisms
	\begin{equation*}
		Z\into X_u\ast Y_v\overset{\ell_{X_u,Y_v}}{\onto}\ell(X_u,Y_v)\isoto p(X_u,Y_v)\into p(X,Y),
	\end{equation*}
	where the first and last maps are inclusions. It follows that $\omega$ is a group homomorphism. We have already observed that the image of $\omega$ is contained in $p_Z(X,Y)$, so we can consider $\omega$ as a homomorphism $Z\to p_Z(X,Y)$. The proof is complete.
\end{proof}

Note that in the above lemma, one has $\omega(g,k)\in\ell_{X_u,Y_v}(g,k)\subseteq \ell_{X,Y}(g,k)$.

\begin{notation}\label{not:omegagammat}
	\textbf{(a)} Set
	\begin{equation*}
		\Omega=\set{\omega\in\Hom\bigl(Z,p_Z(X,Y)\bigr)|\omega(g,k)\in\ell_{X,Y}(g,k)\text{ for all }(g,k)\in Z}.
	\end{equation*}
	In other words, $\Omega$ consists of all group homomorphisms $\omega:Z\to p_Z(X,Y)$ that make the diagram
	\[\begin{tikzcd}
		{p_Z(X,Y)} & {\ell_Z(X,Y)} \\
		Z
		\arrow[two heads, from=1-1, to=1-2]
		\arrow["\omega", from=2-1, to=1-1]
		\arrow["\ell"', two heads, from=2-1, to=1-2]
	\end{tikzcd}\]
	commute.\\
	\textbf{(b)} Let
	\begin{equation*}
		\Gamma=\set{(X',Y')|X'\subgp X,Y'\subgp Y, Z\subgp X'\ast Y'}.
	\end{equation*}
	The set $\Gamma$ is nonempty, since $(X,Y)\in\Gamma$. Note that if $U\in\lset{X}$, $V\in\lset{Y}$, and $u\tensor v\in(U\hat{\tensor}_{X,Y}V)^Z$ then $(X_u,Y_v)\in\Gamma$ by Lemma \ref{lem:stabincomp}.\\
	\textbf{(c)} Set
	\begin{equation*}
		T=\set{(g,h,k)|(g,h)\in X,(h,k)\in Y,(g,k)\in N_{X\ast Y}(Z)}\subgp G\times H\times K.
	\end{equation*}
	The sets $\Omega$ and $\Gamma$ each possess a left action of the group $T$: see \ref{actiont} below.
\end{notation}

\begin{nothing}\label{stroft}\textit{Structure of the group }$T$\textit{.}
	The map $k(X,Y)\into T$, $h\mapsto (1,h,1)$, is an injective group homomorphism whose image is precisely the subgroup of triples in $T$ of the form $(1,h,1)$. For ease, we will identify $k(X,Y)$ with its image under this homomorphism. Note also that we have a short exact sequence
	\begin{equation*}
		1\to k(X,Y)\to T\to N_{X\ast Y}(Z)\to 1\,,
	\end{equation*}
	where the map $T\onto N_{X\ast Y}(Z)$ is defined in the obvious way, mapping $(g,h,k)$ to $(g,k)$. 
\end{nothing}

\begin{nothing}\label{actiont}\textit{The actions of }$T$\textit{ on }$\Omega$\textit{ and }$\Gamma$\textit{.}
	\textbf{(a)} The group $T$ acts on $\Omega$ via
	\begin{equation*}
		{}^{(g,h,k)}\omega=c_h\circ\omega\circ c_{(g,k)}^{-1},
	\end{equation*}
	where $(g,h,k)\in T$ and $\omega\in\Omega$. Indeed, if $(g,h,k)\in T$ then $(g,k)\in N_{X\ast Y}(Z)$ and therefore $\ell_{X,Y}(g,k)\in\ell(X,Y)$ normalizes the image of $Z$ under $\ell_{X,Y}$, which is $\ell_Z(X,Y)=p_Z(X,Y)/k(X,Y)$. Since $h\in\ell_{X,Y}(g,k)$ we have $\ell_{X,Y}(g,k)=h\cdot k(X,Y)$. Hence $h$ normalizes $p_Z(X,Y)$, and we see that if $\omega\in\Omega$ then ${}^{(g,h,k)}\omega$, as defined above, is a homomorphism from $Z$ to $p_Z(X,Y)$. Furthermore, if $(x,z)\in Z$ then we have
	\begin{align*}
		({}^{(g,h,k)}\omega)(x,z)k(X,Y)	&={}^h\omega(x^g,z^k)k(X,Y)={}^{hk(X,Y)}\omega(x^g,z^k)k(X,Y)\\
		&={}^{\ell(g,k)}\ell(x^g,z^k)=\ell(x,z).
	\end{align*}
	Therefore ${}^{(g,h,k)}\omega\in\Omega$.\\
	\textbf{(b)} Also the set $\Gamma$ possesses a left action of $T$, defined by
	\begin{equation*}
		{}^{(g,h,k)}(X',Y')=({}^{(g,h)}X',{}^{(h,k)}Y').
	\end{equation*}
	for each $(g,h,k)\in T$ and $(X',Y')\in\Gamma$. This follows from the identity ${}^{(g,h)}X'\ast {}^{(h,k)}Y'={}^{(g,k)}(X'\ast Y')$, which holds for any $g\in G$, $h\in H$, $k\in K$, and any subgroups $X'\subgp G\times H$, $Y'\subgp H\times K$ (cf. \cite[Lemma 2.2(d)]{Boltje_2020}).
\end{nothing}

\begin{nothing}\label{X(omega)Y(omega)}
	For each $\omega\in\Omega$ set
	\begin{equation*}
		X(\omega)=\set{(g,\omega(g,k))|(g,k)\in Z}\subgp X
	\end{equation*}
	and
	\begin{equation*}
		Y(\omega)=\set{(\omega(g,k),k)|(g,k)\in Z}\subgp Y.
	\end{equation*}
	It is straightforward to verify that $(X(\omega),Y(\omega))\in\Gamma$ for any $\omega\in\Omega$ and that the map $\Omega\to\Gamma$, $\omega\mapsto(X(\omega),Y(\omega))$, is $T$-equivariant.
\end{nothing}

\begin{nothing}\label{N(omega)M(omega)}\textit{The subgroups }$N_{X\ast Y}^{\omega}(Z)$\textit{, }$\tilde{X}(\omega)$\textit{, and }$\tilde{Y}(\omega)$\textit{.}
	\textbf{(a)} For each $\omega\in\Omega$ let $N_{X\ast Y}^{\omega}(Z)$ denote the image of the stabilizer $T_\omega$ in $N_{X\ast Y}(Z)$ under the projection homomorphism in \ref{stroft}. Note that a pair $(g,k)\in N_{X\ast Y}(Z)$ belongs to $N_{X\ast Y}^{\omega}(Z)$ if and only if there exists a linking element $h\in\ell_{X,Y}(g,k)$ such that $c_h\circ\omega\circ c_{(g,k)}^{-1}=\omega$. It is clear then that $Z\subgp N_{X\ast Y}^\omega(Z)$.
	
	The set $k(X,Y)\backslash\Omega$ of orbits of $\Omega$ under the action of the normal subgroup $k(X,Y)\nor T$ inherits an action of $N_{X\ast Y}(Z)\iso T/k(X,Y)$. Note that if $\omega\in\Omega$ then the stabilizer in $N_{X\ast Y}(Z)$ of the $k(X,Y)$-orbit of $\omega$ is precisely $N_{X\ast Y}^{\omega}(Z)$. We record that the $N_{X\ast Y}(Z)$-orbits of $k(X,Y)\backslash\Omega$ are in natural 1-to-1 correspondence with the $T$-orbits of $\Omega$.\\
	\textbf{(b)} If $\omega\in\Omega$ then we set
	\begin{equation*}
		\tilde{X}(\omega)=N_X(X(\omega))\subgp X\qquad\text{and}\qquad \tilde{Y}(\omega)=N_Y(Y(\omega))\subgp Y.
	\end{equation*}

	The following lemma establishes the relationship between the subgroups $N_{X\ast Y}^{\omega}(Z)$, $\tilde{X}(\omega)$, and $\tilde{Y}(\omega)$ just defined.
\end{nothing}

\begin{lemma}\label{containment1}
	If $\omega\in\Omega$ then $N_{X\ast Y}^\omega(Z)\subgp \tilde{X}(\omega)\ast \tilde{Y}(\omega)$.
\end{lemma}

\begin{proof}
	Let $\omega\in\Omega$ and suppose $(g,k)\in N_{X\ast Y}^\omega(Z)$. Then there exists a linking element $h\in\ell_{X,Y}(g,k)$ such that $(g,h,k)\in T$ and ${}^{(g,h,k)}\omega=\omega$. From the $T$-equivariance of the map $\Omega\to\Gamma$ defined in \ref{X(omega)Y(omega)} it follows that 
	\begin{equation*}
		\bigl({}^{(g,h)}X(\omega),{}^{(h,k)}Y(\omega)\bigr)={}^{(g,h,k)}\bigl(X(\omega),Y(\omega)\bigr)=\bigl(X(\omega),Y(\omega)\bigr).
	\end{equation*}
	In other words, $(g,h)\in \tilde{X}(\omega)$ and $(h,k)\in \tilde{Y}(\omega)$. We see then that $(g,k)\in \tilde{X}(\omega)\ast \tilde{Y}(\omega)$. The proof is complete.
\end{proof}

We conclude this section by noting that, with the definitions and notations set above, the statement of Theorem \ref{thm:1} is now precise.

\section{Proof of Theorem \ref{thm:1}}\label{sec:pfthm1}

We continue to let $G$, $H$, and $K$ be finite groups, $X\subgp G\times H$, $Y\subgp H\times K$, and $Z\subgp X\ast Y$. We also let $U\in\lset{X}$ and $V\in\lset{Y}$. The notation of Section \ref{sec:setupthm1} will be in effect throughout this section. 

\begin{nothing}\label{thetaomega}\textit{Definition of the maps }$\theta_\omega$ and $\eta_{U,V}$\textit{.}
	\textbf{(a)} Let $\omega\in\Omega$. Then $U^{X(\omega)}$ is a left $\tilde{X}(\omega)$-set and $V^{Y(\omega)}$ is a left $\tilde{Y}(\omega)$-set, where $X(\omega)$, $Y(\omega)$ are as in \ref{X(omega)Y(omega)} and $\tilde{X}(\omega)$, $\tilde{Y}(\omega)$ are as in \ref{N(omega)M(omega)}(b). By the universal property of the extended tensor product (Proposition \ref{prop:univpropexttens}), there exists a unique $\tilde{X}(\omega)\ast \tilde{Y}(\omega)$-equivariant map
	\begin{equation*}
		U^{X(\omega)}\underset{\tilde{X}(\omega),\tilde{Y}(\omega)}{\hat{\tensor}}V^{Y(\omega)}\to\Res_{\tilde{X}(\omega)\ast \tilde{Y}(\omega)}^{X\ast Y}(U\underset{X,Y}{\hat{\tensor}}V)
	\end{equation*}
	defined by $u\tensor v\mapsto u\tensor v$ for all $u\in U^{X(\omega)}$ and $v\in V^{Y(\omega)}$. Since $N_{X\ast Y}^\omega(Z)\subgp \tilde{X}(\omega)\ast \tilde{Y}(\omega)$ by Lemma \ref{containment1}, this map is also $N_{X\ast Y}^\omega(Z)$-equivariant. Its image is contained in $(U\hat{\tensor}_{X,Y}V)^Z$, since if $u\in U^{X(\omega)}$, $v\in V^{Y(\omega)}$, and $(g,k)\in Z$ then within $U\hat{\tensor}_{X,Y}V$ we have
	\begin{equation*}
		(g,k)\cdot(u\tensor v)=(g,\omega(g,k))u\tensor (\omega(g,k),k)v=u\tensor v.
	\end{equation*}
	Thus, the map above can be considered as a morphism of $N_{X\ast Y}^\omega(Z)$-sets:
	\begin{equation*}
		\zeta_\omega:\Res_{N_{X\ast Y}^\omega(Z)}^{\tilde{X}(\omega)\ast \tilde{Y}(\omega)}
		\Bigl(U^{X(\omega)}\underset{\tilde{X}(\omega),\tilde{Y}(\omega)}{\hat{\tensor}}V^{Y(\omega)}\Bigr)
		\to\Res_{N_{X\ast Y}^\omega(Z)}^{N_{X\ast Y}(Z)}\bigl((U\underset{X,Y}{\hat{\tensor}}V)^Z\bigr).
	\end{equation*}
	By the adjunction of \ref{induction} we obtain an $N_{X\ast Y}(Z)$-equivariant map
	\begin{equation*}
		\theta_\omega=\hat{\zeta_\omega}:\Ind_{N_{X\ast Y}^\omega(Z)}^{N_{X\ast Y}(Z)}\
		\Bigl(\Res_{N_{X\ast Y}^\omega(Z)}^{\tilde{X}(\omega)\ast \tilde{Y}(\omega)}
		\bigl(U^{X(\omega)}\underset{\tilde{X}(\omega),\tilde{Y}(\omega)}{\hat{\tensor}}V^{Y(\omega)}\bigr)\Bigr)
		\to (U\underset{X,Y}{\hat{\tensor}}V)^Z
	\end{equation*}
	defined by
	\begin{equation*}
		\theta_\omega\bigl((g,k)\tensor(u\tensor v)\bigr)=(g,k)\cdot(u\tensor v)
	\end{equation*}
	for all $(g,k)\in N_{X\ast Y}(Z)$, $u\in U^{X(\omega)}$, and $v\in V^{Y(\omega)}$.\\
	\textbf{(b)} Now let $[T\backslash\Omega]$ be a set of representatives of the $T$-orbits of $\Omega$. The coproduct of the maps $\theta_\omega$, $\omega\in[T\backslash\Omega]$, is the map $\eta_{U,V}$ in Theorem~\ref{thm:1}. Note that $\eta_{U,V}$ is defined for all $(U,V)\in\lset{X}\times\lset{Y}$, without restriction on $k(X_u,Y,_v)$ for $(u,v)\in U\times V$. It is straight forward to verify that the collection $\eta$ of maps 
$\eta_{U,V}$, for $(U,V)\in\lset{X}\times \lset{Y}$, is a natural transformation of biadditive functors $\lset{X}\times\lset{Y}\to\lset{N_{X\ast Y}(Z)}$.
\end{nothing}

\begin{lemma}\label{lem:thetainj}
	Let $X\subgp G\times H$, $Y\subgp H\times K$, and $Z\subgp X\ast Y$. Let $U\in\lset{X}$ and $V\in\lset{Y}$ such that $k(X_u,Y_v)=1$ for all $u\in U$ and $v\in V$. Then the map $\theta_\omega$ defined in \ref{thetaomega} is injective for any $\omega\in\Omega$.
\end{lemma}

\begin{proof}
	Fix a homomorphism $\omega\in\Omega$ and let $\zeta_\omega$ be as in \ref{thetaomega}. We want to apply Lemma \ref{lem:injmapfrominduct} to show that $\theta_\omega=\hat{\zeta_\omega}$ is injective. Thus we need to verify both
	\begin{itemize}
		\item[(1)] $\zeta_\omega$ is injective; and
		\item[(2)] Whenever $(g,k)\in N_{X\ast Y}(Z)$, $u,u'\in U^{X(\omega)}$, and $v,v'\in V^{Y(\omega)}$ are such that $(g,k)\zeta_\omega(u\tensor v)=\zeta_\omega(u'\tensor v')$ then $(g,k)\in N_{X\ast Y}^\omega(Z)$.
	\end{itemize}

	We first establish (1). Suppose that $u,u'\in U^{X(\omega)}$ and $v,v'\in V^{Y(\omega)}$ are such that $\zeta_\omega(u\tensor v)=\zeta_\omega(u'\tensor v')$. Then $u\tensor v=u'\tensor v'$ in $U\hat{\tensor}_{X,Y}V$, so by Definition \ref{defn:exttensprod} there exists an element $h\in k(X,Y)$ such that $(1,h)u=u'$ and $(h,1)v=v'$. Now if $(g,k)\in Z$ then $(g,\omega(g,k))\in X(\omega)$ fixes $u$, and therefore ${}^{(1,h)}(g,\omega(g,k))\in X_{u'}$. Since also $(g,\omega(g,k))\in X_{u'}$, we have
	\begin{equation*}
		\bigl(1,{}^h\omega(g,k)\omega(g,k)^{-1}\bigr)={}^{(1,h)}\bigl(g,\omega(g,k)\bigr) \bigl(g,\omega(g,k)\bigr)^{-1}\in X_{u'}.
	\end{equation*}
	Similarly, $\bigl({}^h\omega(g,k)\omega(g,k)^{-1},1\bigr)\in Y_{v'}$. Therefore
	\begin{equation*}
		{}^h\omega(g,k)\omega(g,k)^{-1}\in k(X_{u'},Y_{v'}),
	\end{equation*}
	and since $k(X_{u'},Y_{v'})$ is trivial by assumption it follows that ${}^h\omega(g,k)=\omega(g,k)$. Since this equality holds for any $(g,k)\in Z$ we see that $(1,h)$ centralizes $X(\omega)$ and $(h,1)$ centralizes $Y(\omega)$. In particular, $h\in k(\tilde{X}(\omega),\tilde{Y}(\omega))$. Thus, within the domain $U^{X(\omega)}\hat{\tensor} V^{Y(\omega)}$ of $\zeta_\omega$ we have $u\tensor v=(1,h)u\tensor(h,1)v=u'\tensor v'$. This establishes (1).
	
	It remains to verify (2). Suppose that $(g,k)\in N_{X\ast Y}(Z)$, $u,u'\in U^{X(\omega)}$, and $v,v'\in V^{Y(\omega)}$ are such that $(g,k)\zeta_\omega(u\tensor v)=\zeta_\omega(u'\tensor v')$. Then $(g,k)\cdot(u\tensor v)=u'\tensor v'$ in $U\hat{\tensor}_{X,Y} V$, so there exists a linking element $h\in\ell_{X,Y}(g,k)$ such that $(g,h)u=u'$ and $(h,k)v=v'$. Now if $(x,z)\in Z$ then $(x,z)^{(g,k)}\in Z$, hence $\bigl(x^g,\omega(x^g,z^k)\bigr)\in X(\omega)$ and $\bigl(\omega(x^g,z^k),z^k\bigr)\in Y(\omega)$. Arguing as in the previous paragraph, we have
	\begin{equation*}
		\bigl(1,{}^h\omega(x^g,z^k)\omega(x,z)^{-1}\bigr)={}^{(g,h)}\bigl(x^g,\omega(x^g,z^k)\bigr) \bigl(x,\omega(x,z)\bigr)^{-1}\in X_{u'}
	\end{equation*}
	and
	\begin{equation*}
		\bigl({}^h\omega(x^g,z^k)\omega(x,z)^{-1},1\bigr)={}^{(h,k)}\bigl(\omega(x^g,z^k),z^k\bigr) \bigl(\omega(x,z),z\bigr)^{-1}\in Y_{v'}.
	\end{equation*}
	Since $k(X_{u'},Y_{v'})$ is trivial it follows that ${}^h\omega(x^g,z^k)=\omega(x,z)$. Since $(x,z)$ was an arbitrary element of $Z$ we have $c_h\circ\omega\circ c_{(g,k)}^{-1}=\omega$, i.e., $(g,h,k)\in T_\omega$. We conclude that $(g,k)\in N_{X\ast Y}^\omega(Z)$ and (2) holds. Finally, Lemma \ref{lem:injmapfrominduct} implies that $\theta_\omega=\hat{\zeta_\omega}$ is injective. The proof is complete.
\end{proof}

We are finally in position to establish Theorem \ref{thm:1}.

\begin{proof}[Proof of Theorem \ref{thm:1}.]
	We begin by noting that the map $\eta_{U,V}$ is the coproduct of the $N_{X\ast Y}(Z)$-equivariant maps $\theta_\omega$ (defined in \ref{thetaomega}) as $\omega$ runs over the set $[T\backslash\Omega]$ of representatives for the orbits of $T$ on $\Omega$. So to complete the proof all that must be shown is that $\eta_{U,V}$ is both injective and surjective.
	
	We first show that $\eta_{U,V}$ is injective. Toward this end, suppose that $\omega,\omega'\in\Omega$, $u\in U^{X(\omega)}$, $v\in V^{Y(\omega)}$, $u'\in U^{X(\omega')}$, $v'\in V^{Y(\omega')}$ and $(g,k)\in N_{X\ast Y}(Z)$ are such that $u\tensor v=(g,k)\cdot(u'\tensor v')$ in $U\hat{\tensor}_{X,Y}V$. We claim that $\omega$ and $\omega'$ are $T$-conjugate. Indeed, by Definition \ref{defn:exttensprod} there exists a linking element $h\in\ell_{X,Y}(g,k)$ satisfying $(g,h)u'=u$ and $(h,k)v'=v$. Note that $(g,h,k)\in T$. Now, if $(x,z)\in Z$ then $\bigl({}^gx,\omega({}^gx,{}^kz)\bigr)\in X(\omega)$ fixes $u$ and $\bigl(\omega({}^gx,{}^kz),{}^kz\bigr)\in Y(\omega)$ fixes $v$. Moreover, since $u'\in U^{X(\omega')}$ the element ${}^{(g,h)}\bigl(x,\omega'(x,z)\bigr)$ fixes $u$. Likewise, ${}^{(h,k)}\bigl(\omega'(x,z),z\bigr)$ fixes $v$. It follows that
	\begin{equation*}
		\bigl(1,{}^h\omega'(x,z)\omega({}^gx,{}^kz)^{-1}\bigr) = {}^{(g,h)}\bigl(x,\omega'(x,z)\bigr) \bigl({}^gx,\omega({}^gx,{}^kz)\bigr)^{-1}\in X_u
	\end{equation*}
	and
	\begin{equation*}
		\bigl({}^h\omega'(x,z)\omega({}^gx,{}^kz)^{-1},1\bigr)={}^{(h,k)}\bigl(\omega'(x,z),z\bigr) \bigl(\omega({}^gx,{}^kz),{}^kz\bigr)^{-1}\in Y_v.
	\end{equation*}
	Since $k(X_u,Y_v)=1$ by assumption, we see that ${}^h\omega'(x,z)=\omega({}^gx,{}^kz)$. The element $(x,z)\in Z$ was arbitrary, so $c_h\circ\omega'=\omega\circ c_{(g,k)}$. In other words, ${}^{(g,h,k)}\omega'=\omega$. Thus $\omega$ and $\omega'$ are $T$-conjugate, as claimed.
	
	Now suppose that $\omega,\omega'\in[T\backslash\Omega]$, $(g,k),(g',k')\in N_{X\ast Y}(Z)$, $u\in U^{X(\omega)}$, $v\in V^{Y(\omega)}$, $u'\in U^{X(\omega')}$, and $v'\in V^{Y(\omega')}$ are such that
	\begin{equation*}
		\eta_{U,V}\bigl((g,k)\tensor(u\tensor v)\bigr) = \eta_{U,V}\bigl((g',k')\tensor(u'\tensor v')\bigr).
	\end{equation*}
	Then $\theta_\omega\bigl((g,k)\tensor(u\tensor v)\bigr) = \theta_{\omega'}\bigl((g',k')\tensor(u'\tensor v')\bigr)$. This implies that $(g,k)\cdot(u\tensor v)=(g',k')\cdot(u'\tensor v')$ in $U\hat{\tensor}_{X,Y}V$, so by the claim of the previous paragraph the homomorphisms $\omega$ and $\omega'$ must be $T$-conjugate. But $\omega$ and $\omega'$ belong to $[T\backslash\Omega]$, so in fact $\omega=\omega'$. Now $\theta_\omega$ is injective by Lemma \ref{lem:thetainj}, so we see that $(g,k)\tensor(u\tensor v)=(g',k')\tensor(u'\tensor v')$ within the domain of $\eta_{U,V}$. We have shown that $\eta_{U,V}$ is injective.
	
	Finally, we show that $\eta_{U,V}$ is surjective. Let $u\in U$ and $v\in V$ be such that $u\tensor v\in(U\hat{\tensor}_{X,Y}V)^Z$. Then by Lemma \ref{foretasurj} there exists a homomorphism $\omega\in\Omega$ such that $X(\omega)\subgp X_u$ and $Y(\omega)\subgp Y_v$. Let $\omega'$ be the unique element of $[T\backslash\Omega]$ that is in the $T$-orbit of $\omega$, and let $(g,h,k)\in T$ be such that ${}^{(g,h,k)}\omega=\omega'$. Then ${}^{(g,h)}X(\omega)=X(\omega')$ and ${}^{(h,k)}Y(\omega)=Y(\omega')$ by the $T$-equivariance of the map described in \ref{X(omega)Y(omega)}. It follows that, if we set $u'=(g,h)u$ and $v'=(h,k)v$, then $u'\in U^{X(\omega')}$, $v'\in V^{Y(\omega')}$, and
	\begin{equation*}
		\eta_{U,V}\bigl((g,k)^{-1}\tensor(u'\tensor v')\bigr) =(g,k)^{-1}\cdot\bigl((g,h)u\tensor(h,k)v\bigr)=u\tensor v.
	\end{equation*}
	We have shown that $\eta_{U,V}$ is surjective. The proof is complete.
\end{proof}

\begin{remark}\label{rmk:natural}
If one chooses two sets $[T\backslash\Omega]$ and $[T\backslash\Omega]'$ of representatives of the $T$-orbits of $\Omega$ then one obtains two natural transformations $\eta$ and $\eta'$. We will explain how $\eta$ and $\eta'$ are related. For ease of notation, set
	\begin{equation*}
		W_\omega=\Ind_{N_{X\ast Y}^\omega(Z)}^{N_{X\ast Y}(Z)}
		\Bigl(\Res_{N_{X\ast Y}^\omega(Z)}^{\tilde{X}(\omega)\ast \tilde{Y}(\omega)}
		\bigl(U^{X(\omega)}\underset{\tilde{X}(\omega),\tilde{Y}(\omega)}{\hat{\tensor}}V^{Y(\omega)}\bigr)\Bigr)
	\end{equation*}
	for each $\omega\in\Omega$. If $\omega,\omega'\in\Omega$ belong to the same $T$-orbit then $\omega'=\lexp{(g,h,k)}{\omega}$ for some $(g,h,k)\in T$. It is straightforward to show that the map
	\begin{equation*}
	   \varphi_{(g,h,k)}\colon W_\omega\to W_{\omega'}\,,\quad 
	   (x,z)\tensor(u\tensor v)\mapsto (x,z)(g,k)^{-1}\tensor \bigl((g,h)u\tensor(h,k)v\bigr)\,,
	\end{equation*}
	for any $(x,z)\in N_{X\ast Y}(Z)$, $u\in U^{X(\omega)}$, and $v\in V^{Y(\omega)}$, is a well-defined isomorphism of $N_{X\ast Y}(Z)$-sets. It is also straightforward to show that $\varphi_{(g,h,k)}$ does not depend on the element $(g,h,k)\in T$ with $\omega'=\lexp{(g,h,k)}{\omega}$. Thus, we may denote this map by $\varphi_{\omega',\omega}$. Then we have $\theta_{\omega'}\circ\varphi_{\omega',\omega}=\theta_{\omega} \colon W_\omega\to (U \hat{\tensor}_{X,Y}V)^Z$. Therefore, if $\omega$ runs through $[T\backslash\Omega]$ and $\omega'$ denotes the unique element in $[T\backslash\Omega]'$ which is in the same $T$-orbit as $\omega$, then the coproduct of the maps $\varphi_{\omega',\omega}$ defines a natural isomorphism between the domains of $\eta_{U,V}$ and $\eta'_{U,V}$ such that $\eta'\circ \varphi = \eta$.
\end{remark}

\section{Trivial Source Modules and Brauer Pairs}\label{sec:trivsourcemodsandbrauerpairs}

For the remainder of this article $F$ denotes a field of positive characteristic $p$. We may sometimes need to assume that $F$ is ``large enough'' for the finite group $G$, which means that $F$ contains a primitive root of unity whose order is equal to $\exp(G)_{p'}$, the $p'$-part of the exponent of $G$. In this section we review various definitions and facts relating to the modular representation theory of finite groups, in preparation for the two consequences of Theorem \ref{thm:1} we present in Section \ref{sec:apps}.

\begin{nothing}\label{FGmods}$FG$\textit{-modules.}
	\textbf{(a)} The categories of finitely-generated $FG$-modules and $(FG,FH)$-bimodules are denoted by $\lmod{FG}$ and $\bimod{FG}{FH}$, respectively. As in the case of bisets, any $(FG,FH)$-bimodule $M$ can be considered a left $F[G\times H]$-module (and vice versa) by defining $(g,h)m=gmh^{-1}$ for any $(g,h)\in G\times H$ and $m\in M$, and ``switching sides'' in this way defines an isomorphism of categories $\bimod{FG}{FH}\iso\lmod{F[G\times H]}$. In the sequel, all modules and bimodules for finite group algebras will be assumed finitely-generated.\\
	\textbf{(b)} If $M\in\lmod{FG}$ and $H\subgp G$ then $M^H$ denotes the $FN_G(H)$-submodule of $H$-fixed points of $M$.\\
	\textbf{(c)} Recall that we have a ``linearization'' functor $\lset{G}\to\lmod{FG}$ mapping a finite left $G$-set $U$ to the permutation $FG$-module $FU$, which has an $F$-basis formed by the elements of $U$.
\end{nothing}

\begin{nothing}\textit{Blocks.} (cf. \cite[Chapter 6]{Linckelmann_2018_2})
	Let $G$ be a finite group. Recall that the primitive idempotents of the center $Z(FG)$ of $FG$ are called the \textit{block idempotents} of the group algebra $FG$. If $e$ is a block idempotent of $FG$ then the corresponding \textit{block algebra} is the ideal $FGe$ generated by $e$, which is in fact an $F$-algebra with identity $e$. If $B$ is a block algebra of $FG$ then we may denote the identity (block idempotent) of $B$ by $e_B$. The set of block idempotents of $FG$ will be denoted by $\bli(FG)$ and the set of block algebras of $FG$ will be denoted by $\Bl(FG)$. In the sequel the term ``block'' may refer to either a block algebra or a block idempotent. 
	
	The antipode map 
	\begin{equation*}
		-^\ast:FG\to FG,\qquad \sum_{g\in G}\alpha_gg\mapsto\sum_{g\in G}\alpha_gg^{-1},
	\end{equation*}
	induces bijections on the sets of block idempotents and block algebras. If $H$ is another finite group and $F$ is large enough for $G\times H$ then every block idempotent of $F[G\times H]\iso FG\tensor_F FH$ is of the form $e\tensor f^\ast$ for uniquely determined block idempotents $e\in\bli(FG)$, $f\in\bli(FH)$. In particular, every block algebra of $F[G\times H]$ is of the form $A\tensor B^\ast$ for uniquely determined block algebras $A\in\Bl(FG)$, $B\in\Bl(FH)$ (note that we may suppress the subscript in the notation $\tensor_F$ when discussing tensor products of blocks).
	
	Recall that a block $B$ of $FG$ determines a full conjugacy class of $p$-subgroups of $G$ called the \textit{defect groups} of $B$ (cf. \cite[Definition 6.1.1]{Linckelmann_2018_2}). Recall also that an $FG$-module $M$ is said to \textit{belong} to the block $B$ if $e_BM=M$, or equivalently, if $e_Bm=m$ for all $m\in M$. More generally, if $e$ is a central idempotent of $FG$ then $M$ \textit{belongs to} $e$ if $em=m$ for all $m\in M$. 
\end{nothing}

\begin{nothing}\textit{Trivial source modules.} (cf. \cite{Broue_1985}, \cite[Chapter 5.10]{Linckelmann_2018})
	Recall that a module $M\in\lmod{FG}$ is a \textit{trivial source} or $p$\textit{-permutation} module if every indecomposable direct summand of $M$ has the trivial module $F$ as a source or, equivalently, if $M$ is isomorphic to a direct summand of a permutation $FG$-module. A finite-dimensional $(FG,FH)$-bimodule $M$ is a \textit{trivial source bimodule} if $M$ is trivial source when viewed as a left $F[G\times H]$-module. The full subcategory of $\lmod{FG}$ (respectively, of $\bimod{FG}{FH}$) whose objects are the trivial source $FG$-modules (resp., trivial source $(FG,FH)$-bimodules) is denoted by $\ltriv{FG}$ (resp., $\bitriv{FG}{FH}$). The category $\ltriv{FG}$ is the idempotent completion of $\lperm{FG}$, the category of finite-dimensional permutation $FG$-modules, and is an $F$-linear, symmetric monoidal category with respect to the tensor product $-\tensor_F-$. The class of trivial source modules is closed under the operations of induction, restriction, inflation, and conjugation. Moreover, if $M$ is a trivial source $FG$-module then the dual $M^\circ=\Hom_F(M,F)$ is also a trivial source $FG$-module, with the action of an element $g\in G$ on $\varphi\in M^\circ$ defined by $(g\varphi)(m)=\varphi(g^{-1}m)$ for all $m\in M$. More generally, if $M$ is a trivial source $(FG,FH)$-bimodule then the dual $M^\circ$ is a trivial source $(FH,FG)$-bimodule via $(h\varphi g)(m)=\varphi(gmh)$ for all $h\in H$, $g\in G$, $\varphi\in M^\circ$ and $m\in M$.
	
	If $e$ is a central idempotent of $FG$ then a \textit{trivial source} $FGe$-module is simply a trivial source $FG$-module that belongs to $e$. The full subcategory of $\ltriv{FG}$ whose objects are the trivial source $FGe$-modules is denoted by $\ltriv{FGe}$. This subcategory is closed under direct sums and taking direct summands.
\end{nothing}

\begin{nothing}\textit{The trivial source ring.} (cf. \cite[Chapter 5.5]{Benson_1991_vol1})
	The Grothendieck ring of $\ltriv{FG}$ (with respect to split short exact sequences) is called the \textit{trivial source ring} and is denoted by $T_F(G)=T(FG)$. As an abelian group, $T_F(G)$ is free of finite rank with a \textit{standard basis} consisting of the isomorphism classes $[M]$ of indecomposable trivial source $FG$-modules $M$. Each of the operations on trivial source modules mentioned in the previous paragraph induce maps between trivial source rings. For example, one has a unital ring automorphism $(-)^\circ:T_F(G)\isoto T_F(G)$ mapping the isomorphism class $[M]$ of a trivial source $FG$-module $M$ to $[M^\circ]$.
	
	The Grothendieck group of $\bitriv{FG}{FH}$ is denoted by $T_F(G,H)=T(FG,FH)$. As ``switching sides'' defines an isomorphism $\bitriv{FG}{FH}\iso\ltriv{F[G\times H]}$, the group $T_F(G,H)$ is isomorphic to $T_F(G\times H)$. Via this isomorphism we sometimes view the elements of $T_F(G,H)$ as virtual trivial source $F[G\times H]$-modules (rather than virtual trivial source bimodules). Note that ``taking duals'' defines a homomorphism $(-)^\circ:T_F(G,H)\to T_F(H,G)$. 
	
	If $e$ is a central idempotent of $FG$ then the Grothendieck group of $\ltriv{FGe}$ is denoted by $T_F(G,e)=T(FGe)$. Note that $T_F(G)=\dsum_{B\in\Bl(FG)}T_F(B)$. If $\mu\in T_F(G)$ and $B$ is a block algebra of $FG$ with identity $e_B$ then we write $e_B\mu$ for the image of $\mu$ under the canonical projection $T_F(G)\onto T_F(B)$.
	
	If $H$ is another finite group and $f$ is a central idempotent of $FH$ then the Grothendieck group of the category $\bitriv{FGe}{FHf}$ of trivial source $(FGe,FHf)$-bimodules is denoted by $T(FGe,FHf)$. Once again, ``switching sides'' induces an isomorphism $T(FGe,FHf)\iso T_F(G\times H,e\tensor f^\ast)$. We often identify these groups via this isomorphism without comment. The subgroup of $T(FGe,FHf)$ spanned by the isomorphism classes of indecomposable trivial source $F[G\times H](e\tensor f^\ast)$-modules that have twisted diagonal vertices is denoted by $T^\Delta(FGe,FHf)$. 
\end{nothing}

\begin{nothing}\label{Brauerconst}\textit{The Brauer construction.} (cf. \cite{Broue_1985}, \cite[Chapter 5.4]{Linckelmann_2018})
	If $M\in\lmod{FG}$ and $P$ is a $p$-subgroup of $G$ then the \textit{Brauer construction} of $M$ at $P$ is the $FN_G(P)$-module
	\begin{equation*}
		M(P)=M^P/\sum_{Q\lneq P}\tr_Q^P(M^Q),
	\end{equation*}
	where, for each proper subgroup $Q$ of $P$, $\tr_Q^P:M^Q\to M^P$ denotes the relative trace map (cf. \cite[Chapter 2.5]{Linckelmann_2018}). The canonical surjection $M^P\onto M(P)$ is denoted by $\Br_P^M$. If $A$ is a $G$-algebra over $F$ then $A(P)$ inherits an $N_G(P)$-algebra structure. The Brauer construction at $P$ defines a lax monoidal $F$-linear functor $-(P):\lmod{FG}\to\lmod{FN_G(P)}$. 
\end{nothing}

\begin{nothing}\label{Brauerconsttriv}\textit{The Brauer construction applied to trivial source modules.} (cf. \cite[Chapters 5.8, 5.10]{Linckelmann_2018})
	If $M=FU$ is a permutation $FG$-module, where $U\in\lset{G}$, and $P$ is a $p$-subgroup of $G$ then the composite of the maps $F[U^P]\into (FU)^P\onto (FU)(P)$ is an isomorphism of $FN_G(P)$-modules. In particular, $M(P)$ is a permutation $FN_G(P)$-module. It follows that if $M$ is a trivial source $FG$-module then $M(P)$ is a trivial source $FN_G(P)$-module. The restriction of the functor $-(P)$ to $\ltriv{FG}$ is a strong monoidal functor $\ltriv{FG}\to\ltriv{FN_G(P)}$. Furthermore one has a commutative (up to a natural isomorphism) diagram of functors
	\[\begin{tikzcd}
		{\lset{G}} & {\lset{N_G(P)}} \\
		{\ltriv{FG}} & {\ltriv{FN_G(P)}}
		\arrow["{(-)^P}", from=1-1, to=1-2]
		\arrow[from=1-1, to=2-1]
		\arrow[from=1-2, to=2-2]
		\arrow["{-(P)}"', from=2-1, to=2-2]
	\end{tikzcd}\]
	Note that the Brauer construction induces a unital ring homomorphism $-(P):T_F(G)\to T_F(N_G(P))$.
\end{nothing}

\begin{lemma}\label{vtxbrauer}(\cite[Proposition 5.10.3]{Linckelmann_2018})
	Let $P$ be a $p$-subgroup of $G$, let $M$ be an indecomposable trivial source $FG$-module, and let $Q$ be a vertex of $M$. Then $M(P)\neq\set{0}$ if and only if $P$ is $G$-conjugate to a subgroup of $Q$.
\end{lemma}

\begin{nothing}\textit{The Brauer homomorphism.} (cf. \cite[Theorem 5.4.1]{Linckelmann_2018})
	Let $P$ be a $p$-subgroup of $G$. The $F$-linear projection 
	\begin{equation*}
		FG\onto FC_G(P),\qquad \sum_{g\in G}\alpha_g g\mapsto\sum_{g\in C_G(P)}\alpha_g g,
	\end{equation*}
	restricts to a split surjective homomorphism of $N_G(P)$-algebras
	\begin{equation*}
		\br_P^G:(FG)^P\onto FC_G(P)
	\end{equation*}
	called the \textit{Brauer homomorphism}. Here, $(FG)^P$ is the algebra of $P$-fixed points of $FG$ for the action of $P$ on $FG$ by conjugation. When the group $G$ is contextually clear we may write $\br_P$ in place of $\br_P^G$.
	
	The Brauer homomorphism can be viewed as a special case of the Brauer construction: considering the group algebra $FG$ as a $G$-algebra via conjugation, one has $G^P=C_G(P)$ and $(FG)(P)\iso FC_G(P)$ by the remarks in \ref{Brauerconsttriv}. 
	
	Note that if $e$ is a central idempotent of $FG$ then $\br_P(e)$ is both an $N_G(P)$-stable central idempotent of $FC_G(P)$ and a central idempotent of $FN_G(P)$.
\end{nothing}

The following lemma is well known.

\begin{lemma}\label{M(P)belongsbr(e)}
	Let $e$ be a central idempotent of $FG$ and let $P$ be a $p$-subgroup of $G$. If $M$ is a trivial source $FGe$-module then $M(P)$ belongs to the central idempotent $\br_P(e)$ of $FN_G(P)$.
\end{lemma}

For the remainder of this section we assume that $F$ is ``large enough'' for all of the finite groups considered.

\begin{nothing}\textit{Brauer pairs.} (cf. \cite{Alperin_1979}, \cite[Part IV, Chapter 2]{Aschbacher_2011})
	Recall that an $FG$-\textit{Brauer pair} is a pair $(P,e)$ consisting of a $p$-subgroup $P\subgp G$ and a block idempotent $e$ of $FC_G(P)$. The group $G$ acts by simultaneous conjugation on the set $\BP(FG)$ of $FG$-Brauer pairs. If $(P,e)\in\BP(FG)$ then the stabilizer of $(P,e)$ in $G$ is denoted by $N_G(P,e)$. Note that we always have $PC_G(P)\subgp N_G(P,e)\subgp N_G(P)$. Note also that $e$ is a block idempotent of $FN_G(P,e)$.
	
	The set of $FG$-Brauer pairs is partially ordered: if $(Q,f),(P,e)\in\BP(FG)$ write $(Q,f)\nor(P,e)$ if $Q\subgp P\subgp N_G(Q,f)$ and $\br_P^G(f)e\neq 0$. Then the transitive closure of this relation defines a partial order $\leq$ on $\BP(FG)$ which is compatible with the conjugation action of $G$ (in other words, $\BP(FG)$ is a $G$-poset). If $(P,e)\in\BP(FG)$ and $Q\subgp P$ then there exists a unique block idempotent $f$ of $FC_G(Q)$ such that $(Q,f)\leq (P,e)$ \cite[Theorem 3.4]{Alperin_1979}.
	
	Let $B=FGe_B$ be a block of $FG$. A Brauer pair $(P,e)$ of $FG$ is a $B$\textit{-Brauer pair} if $\br_P^G(e_B)e=e$, or equivalently, if $\br_P^G(e_B)e\neq 0$. The set $\BP(B)$ of $B$-Brauer pairs is $G$-stable and if $(Q,f)\leq(P,e)$ is a containment of Brauer pairs then $(Q,f)$ is a $B$-Brauer pair if and only if $(P,e)$ is a $B$-Brauer pair. The maximal $B$-Brauer pairs form a full $G$-conjugacy class, and a $B$-Brauer pair $(D,e)$ is maximal if and only if $D$ is a defect group of $B$. Note that if $(P,e)\in\BP(B)$ then $(P,e^\ast)\in\BP(B^\ast)$.
\end{nothing}

\begin{nothing}\label{FongReynolds}\textit{The Fong-Reynolds correspondence.}
	Let $(P,e)$ be an $FG$-Brauer pair and set $I=N_G(P,e)$. The sum of the distinct $N_G(P)$-conjugates of $e$ is a block idempotent $\hat{e}$ of $FN_G(P)$. In fact, $\hat{e}$ is the unique block idempotent of $FN_G(P)$ that \textit{covers} $e$, i.e., satisfies $e\cdot\hat{e}\neq 0$. Since $e\in FN_G(P)\hat{e}$ and the ideal of $FN_G(P)\hat{e}$ generated by $e$ is equal to $FN_G(P)\hat{e}$, the block algebras $FN_G(P)\hat{e}$ and $FIe$ are Morita equivalent (see, e.g., \cite[Theorem 2.8.7]{Linckelmann_2018}). In fact, induction defines an equivalence $\Ind_I^{N_G(P)}:\lmod{FIe}\to\lmod{FN_G(P)\hat{e}}$ with inverse $e\Res_I^{N_G(P)}$. Note that these functors restrict to an equivalence between $\ltriv{FN_G(P)\hat{e}}$ and $\ltriv{FIe}$. It follows that $e\Res_I^{N_G(P)}:T(FN_G(P)\hat{e})\isoto T(FIe)$ is a group isomorphism with inverse $\Ind_I^{N_G(P)}$.
\end{nothing}

The following two lemmas will be of particular use in the proof of Theorem \ref{thm:3}.

\begin{lemma}\label{lem:forthm3}
	Let $(P,e)$ be an $FG$-Brauer pair, let $I=N_G(P,e)$, and let $\hat{e}$ be the unique block idempotent of $FN_G(P)$ that covers $e$ (as in \ref{FongReynolds}). Suppose that $H$ is a subgroup of $I$ containing $PC_G(P)$. Then $e$ is a block idempotent of $FH$. Moreover, if $M$ is an $FHe$-module then $\Ind_H^{N_G(P)}(M)$ belongs to $\hat{e}$.
\end{lemma}

\begin{proof}
	Clearly, $e$ is a central idempotent of $FH$. Since $P$ is a normal $p$-subgroup of $H$, \cite[Theorem 5.2.8(ii)]{Nagao_1989} implies that $e$ is a block idempotent of $FH$. Let $M$ be an $FHe$-module. Since $\hat{e}\in FC_G(P)$, to complete the proof it suffices to show that $\hat{e}$ acts trivially on $\Res_{C_G(P)}^{N_G(P)}(\Ind_H^{N_G(P)}(M))$. Now,
	\begin{equation*}
		\Res_{C_G(P)}^{N_G(P)}(\Ind_H^{N_G(P)}(M))\iso\bigoplus_{x\in [N_G(P)/H]}{}^x\Res_{C_G(P)}^H(M)
	\end{equation*}
	by an easy application of the Mackey formula. If $x\in N_G(P)$ then ${}^x\Res_{C_G(P)}^H(M)$ belongs to ${}^xe$. Since ${}^xe\cdot\hat{e}\neq 0$, we see that $\hat{e}$ acts trivially on ${}^x\Res_{C_G(P)}^H(M)$. It follows from the decomposition above that $\hat{e}$ acts trivially on $\Res_{C_G(P)}^{N_G(P)}(\Ind_H^{N_G(P)}(M))$. The proof is complete.
\end{proof}

\begin{lemma}\label{lem:FRcor}
	Let $P$ be a $p$-subgroup of $G$, let $f$ be a central idempotent of $FN_G(P)$, and let $\mathcal{S}$ denote a set of representatives for the $N_G(P)$-conjugacy classes of block idempotents $e\in\bli(FC_G(P))$ such that $ef\neq 0$. Then for any $\mu\in T(FN_G(P)f)$ we have
	\begin{equation*}
		\mu=\sum_{e\in\mathcal{S}}\Ind_{N_G(P,e)}^{N_G(P)}(e\Res_{N_G(P,e)}^{N_G(P)}(\mu)).
	\end{equation*}
\end{lemma}

\begin{proof}
	For each $e\in\mathcal{S}$ let $\hat{e}$ denote the unique block idempotent of $FN_G(P)$ that covers $e$. Then $f=\sum_{e\in\mathcal{S}}\hat{e}$ by \cite[Theorem 5.2.8(ii)]{Nagao_1989}. By the Fong-Reynolds correspondence (\ref{FongReynolds}) the maps $e\Res_{N_G(P,e)}^{N_G(P)}$, $e\in\mathcal{S}$, form the components of an isomorphism
	\begin{equation*}
		T(FN_G(P)f)=\bigoplus_{e\in\mathcal{S}}T(FN_G(P)\hat{e})\isoto\bigoplus_{e\in\mathcal{S}}T(FN_G(P,e)e),
	\end{equation*}
	with inverse $\dsum_{e\in\mathcal{S}}\Ind_{N_G(P,e)}^{N_G(P)}$. Therefore, for any $\mu\in T(FN_G(P)f)$ we have
	\begin{equation*}
		\mu=\sum_{e\in\mathcal{S}}\Ind_{N_G(P,e)}^{N_G(P)}(e\Res_{N_G(P,e)}^{N_G(P)}(\hat{e}\mu)).
	\end{equation*}
	The result follows after noting that $e\Res_{N_G(P,e)}^{N_G(P)}(\hat{e}\mu)=e\Res_{N_G(P,e)}^{N_G(P)}(\mu)$.
\end{proof}

\begin{nothing}\label{ppermequivs}$p$\textit{-Permutation equivalences.} (\cite[Definition 9.8]{Boltje_2020})
	Let $e$ and $f$ be nonzero central idempotents of $FG$ and $FH$. Then a $p$\textit{-permutation equivalence} between $FGe$ and $FHf$ is an element $\gamma\in T^\Delta(FGe,FHf)$ that satisfies $\gamma\tensor_{FH}\gamma^\circ=[FGe]\in T(FGe,FGE)$ and $\gamma^\circ\tensor_{FG}\gamma=[FHf]\in T(FHf,FHf)$.
\end{nothing}

\begin{nothing}\label{gammabrauerpairs}$\gamma$\textit{-Brauer pairs.}
	\textbf{(a)} Let $(P,e)$ be an $FG$-Brauer pair and set $I=N_G(P,e)$. If $M$ is a trivial source $FG$-module then $M(P,e)$ denotes the trivial source $FIe$-module $e\Res_I^{N_G(P)}M(P)$. One obtains a functor $-(P,e):\ltriv{FG}\to\ltriv{FIe}$ and a group homomorphism $-(P,e):T(FG)\to T(FIe)$. If $\mu\in T(FG)$ and $\mu(P,e)\neq 0$, Boltje and Perepelitsky call $(P,e)$ a $\mu$\textit{-Brauer pair} \cite[Definition 9.6]{Boltje_2020}.\\
	\textbf{(b)} Now let $A=FGe_A$ and $B=FHf_B$ be blocks of $FG$ and $FH$, respectively, and suppose that $\gamma\in T^\Delta(A,B)$ is a $p$-permutation equivalence between $A$ and $B$. Then by \cite[Theorem 10.11]{Boltje_2020} the set of $\gamma$-Brauer pairs is a $G\times H$-stable ideal in the partially ordered set of $A\tensor B^\ast$-Brauer pairs, and any two maximal $\gamma$-Brauer pairs are $G\times H$-conjugate. Moreover, any maximal $\gamma$-Brauer pair is of the form $(\Delta(D,\phi,E),e_D\tensor f_E^\ast)$ where $(D,e_D)$ is a maximal $A$-Brauer pair, $(E,f_E)$ is a maximal $B$-Brauer pair, and $\phi:E\isoto D$ is a group isomorphism that induces an isomorphism between the corresponding block fusion systems of $A$ and $B$ (in particular, $\Delta(D,\phi,E)$ is a twisted diagonal subgroup of $G\times H$ --- see \ref{twisteddiagsubg} above). Note that the Brauer subpairs of $(\Delta(D,\phi,E),e_D\tensor f_E^\ast)$ are each uniquely of the form $(\Delta(P,\phi,Q),e_P\tensor f_Q^\ast)$ where $Q\subgp E$, $P=\phi(Q)$, and $e_P,f_Q$ denote the unique block idempotents of $FC_G(P)$ and $FC_H(Q)$, respectively, satisfying $(P,e_P)\leq(D,e_D)$ and $(Q,f_Q)\leq(E,f_E)$.
\end{nothing}

\section{Applications}\label{sec:apps}

With the background definitions and results of Section \ref{sec:trivsourcemodsandbrauerpairs} in hand, we now demonstrate two consequences of Theorem \ref{thm:1}. The first of these results (Theorem \ref{thm:2}) is essentially the ``linearized'' version of Theorem \ref{thm:1} and thus incorporates the extended tensor product of modules, rather than sets. We begin by recalling how this operation is defined. Continue to let $F$ denote a field of positive characteristic $p$.

\begin{nothing}\label{exttensmods}\textit{The Extended Tensor Product of Modules.} (cf. \cite[Section 6]{Boltje_2020})
	\textbf{(a)} Let $X\subgp G\times H$ and $Y\subgp H\times K$. If $M\in\lmod{FX}$ and $N\in\lmod{FY}$ then the \textit{extended tensor product} of $M$ and $N$ is the $F[X\ast Y]$-module $M\hat{\tensor}_{X,Y}N$ whose underlying $F$-vector space is $M\tensor_{F[k(X,Y)]}N$, and with the action of $X\ast Y$ defined by
	\begin{equation*}
		(g,k)\cdot(m\tensor n)=(g,h)m\tensor(h,k)n
	\end{equation*}
	for all $(g,k)\in X\ast Y$, $m\in M$ and $n\in N$, where $h\in H$ is any element linking $g$ and $k$ as in \ref{linkelts}. The extended tensor product is characterized by the universal property of Proposition \ref{prop:univpropexttensmods} below, and defines a functor $\lmod{FX}\times\lmod{FY}\to\lmod{F[X\ast Y]}$ that respects direct sums.\\
	\textbf{(b)} As in \ref{rmksondefnext}(c), there is a canonical isomorphism
	\begin{equation*}
		\Res_{k_1(X)\times k_2(Y)}^{X\ast Y}(M\underset{X,Y}{\hat{\tensor}}N)\iso(\Res_{k_1(X)\times k(X,Y)}^X M)\underset{F[k(X,Y)]}{\tensor}(\Res_{k(X,Y)\times k_2(Y)}^Y N)
	\end{equation*}
	for any $M\in\lmod{FX}$ and any $N\in\lmod{FY}$. In particular, if $X=G\times H$ and $Y=H\times K$ then we have
	\begin{equation*}
		M\underset{G\times H,H\times K}{\hat{\tensor}}N\iso M\underset{FH}{\tensor}N
	\end{equation*}
	for any $(FG,FH)$-bimodule $M$ and any $(FH,FK)$-bimodule $N$. (Recall from \ref{FGmods}(a) that we may view any $(FG,FH)$-bimodule as a left $F[G\times H]$-module after ``switching sides.'')\\
	\textbf{(c)} It is clear that if $U\in\lset{X}$ and $V\in\lset{Y}$ then $$(FU)\underset{X,Y}{\hat{\tensor}} (FV)\iso F[U\underset{X,Y}{\hat{\tensor}} V].$$ In fact, the diagram of functors
	\[\begin{tikzcd}
		{\lset{X}\times\lset{Y}} & {\lset{X\ast Y}} \\
		{\lmod{FX}\times\lmod{FY}} & {\lmod{F[X\ast Y]}}
		\arrow["{-\underset{X,Y}{\hat{\tensor}}-}", from=1-1, to=1-2]
		\arrow[from=1-1, to=2-1]
		\arrow[from=1-2, to=2-2]
		\arrow["{-\underset{X,Y}{\hat{\tensor}}-}", from=2-1, to=2-2]
	\end{tikzcd}\]
	is commutative, up to a natural isomorphism. It follows that if $M\in\ltriv{FX}$ and $N\in\ltriv{FY}$ then $M\hat{\tensor}_{X,Y} N\in\ltriv{F[X\ast Y]}$ (see also \cite[Lemma 7.2(a)]{Boltje_2020}). Moreover, the extended tensor product induces a biadditive map
	\begin{equation*}
		-\underset{X,Y}{\hat{\tensor}}-:T_F(X)\times T_F(Y)\to T_F(X\ast Y).
	\end{equation*}
\end{nothing}

\begin{proposition}\label{prop:univpropexttensmods}
	Let $X\subgp G\times H$, $Y\subgp H\times K$, $M\in\lmod{FX}$, and $N\in\lmod{FY}$. If $L\in\lmod{F[X\ast Y]}$ and $\varphi:M\times N\to L$ is an $F$-bilinear map satisfying
	\begin{equation*}
		(g,k)\varphi(m,n)=\varphi((g,h)m,(h,k)n)
	\end{equation*}
	for all $(g,k)\in X\ast Y$, $h\in\ell_{X,Y}(g,k)$ and $(m,n)\in M\times N$, then there exists a unique $F[X\ast Y]$-homomorphism $\overline{\varphi}:M\hat{\tensor}_{X,Y}N\to L$ making the diagram
	\[\begin{tikzcd}
		{M\times N} & L \\
		{M\underset{X,Y}{\hat{\tensor}}N}
		\arrow["\varphi", from=1-1, to=1-2]
		\arrow["{-\tensor-}"', from=1-1, to=2-1]
		\arrow["{\overline{\varphi}}"', dashed, from=2-1, to=1-2]
	\end{tikzcd}\]
	commute. For $(m,n)\in M\times N$ one has $\overline{\varphi}(m\otimes n) = \varphi(m,n)$.
\end{proposition}

\begin{proof}
	This follows from an argument completely analogous to the one employed in the proof of Proposition \ref{prop:univpropexttens}.
\end{proof}

\begin{lemma}\label{lem:forthm2}
	Let $X\subgp G\times H$, $Y\subgp H\times K$. Let $M$ be an indecomposable trivial source $FX$-module and let $N$ be an indecomposable trivial source $FY$-module. Then the following are equivalent:
	\begin{itemize}
		\item[(1)] $k(P,Q)=1$ for any vertex $P$ of $M$ and any vertex $Q$ of $N$.
		\item[(2)] There exists a finite left $X$-set $U$ and a finite left $Y$-set $V$ such that $M$ is isomorphic to a direct summand of $FU$, $N$ is isomorphic to a direct summand of $FV$, and $k(X_u,Y_v)=1$ for all $u\in U$ and $v\in V$.
	\end{itemize}
\end{lemma}

\begin{proof}
	First assume (1) holds. Let $P$ be a vertex of $M$, $Q$ a vertex of $N$. Set $U=X/P$ and $V=Y/Q$. Since $P$ is a vertex of the trivial source module $M$, the module $\Ind_P^X(F)\iso FU$ has a direct summand isomorphic to $M$. Likewise $FV$ has a direct summand isomorphic to $N$. If $u\in U$ then $X_u$ is conjugate to $P$; hence, $X_u$ is a vertex of $M$. Likewise, if $v\in V$ then $Y_v$ is a vertex of $N$. Therefore $k(X_u,Y_v)=1$ by assumption, and (2) holds.
	
	Now assume (2) and let $U\in\lset{X}$, $V\in\lset{Y}$ be as in the condition. If $P$ is a vertex of $M$ then $M(P)\neq 0$ by Lemma \ref{vtxbrauer}. Since $M$ is isomorphic to a direct summand of $FU$ it follows that $(FU)(P)\neq 0$. As noted in \ref{Brauerconsttriv} we have $(FU)(P)\iso F[U^P]$, so $U^P$ is nonempty. Similarly, if $Q$ is a vertex of $N$ then $V^Q$ is nonempty. Let $u\in U^P$ and $v\in V^Q$. Then $P\subgp X_u$, $Q\subgp Y_v$, and $k(P,Q)\subgp k(X_u,Y_v)=1$. It follows that $k(P,Q)=1$ and condition (1) holds.
\end{proof}

We now present our first application of Theorem \ref{thm:1}. In what follows, we make use of the notations and definitions of Section \ref{sec:setupthm1}.

\begin{theorem}\label{thm:2}
	Let $X\subgp G\times H$, $Y\subgp H\times K$ and let $Z$ be a $p$-subgroup of $X\ast Y$. Let $M\in\ltriv{FX}$, $N\in\ltriv{FY}$ be such that $k(P,Q)=1$ whenever $P$ is a vertex of an indecomposable direct summand of $M$ and $Q$ is a vertex of an indecomposable direct summand of $N$. Let $[T\backslash\Omega]$ denote a set of representatives for the orbits of $T$ on $\Omega$. Then there is a natural isomorphism of $FN_{X\ast Y}(Z)$-modules
	\begin{equation*}
		\eta_{M,N}\colon\bigoplus_{w\in[T\backslash\Omega]}\Ind_{N_{X\ast Y}^{\omega}(Z)}^{N_{X\ast Y}(Z)}(\Res_{N_{X\ast Y}^{\omega}(Z)}^{\tilde{X}(\omega)\ast\tilde{Y}(\omega)}(M(X(\omega))\underset{\tilde{X}(\omega),\tilde{Y}(\omega)}{\hat{\tensor}} N(Y(\omega))))\isoto(M\underset{X,Y}{\hat{\tensor}} N)(Z)
	\end{equation*}
	defined by
	\begin{equation*}
		\eta_{M,N}((g,k)\tensor(\Br_{X(\omega)}^M(m)\tensor\Br_{Y(\omega)}^N(n)))=(g,k)\cdot\Br_Z^{M\hat{\tensor} N}(m\tensor n)
	\end{equation*}
	for all $\omega\in[T\backslash\Omega]$, $(g,k)\in N_{X\ast Y}(Z)$, $m\in M^{X(\omega)}$, and $n\in N^{Y(\omega)}$.
\end{theorem}

\begin{proof}
	We first show that $\eta_{M,N}$ is a well-defined homomorphism between well-defined $FN_{X\ast Y}(Z)$-modules. Let $\omega\in\Omega$. Then $X(\omega)$ and $Y(\omega)$ are $p$-groups, since both are homomorphic images of $Z$. Therefore $M(X(\omega))$, the Brauer construction of $M$ at $X(\omega)$, is a well-defined $F\tilde{X}(\omega)$-module (see \ref{Brauerconst} above). Likewise, $N(Y(\omega))$ is an $F\tilde{Y}(\omega)$-module, and, hence, the extended tensor product of $M(X(\omega))$ and $N(Y(\omega))$ is an $F[\tilde{X}(\omega)\ast\tilde{Y}(\omega)]$-module. By Lemma \ref{containment1} we have $N_{X\ast Y}^\omega(Z)\subgp\tilde{X}(\omega)\ast\tilde{Y}(\omega)$, and it is clear that $N_{X\ast Y}^\omega(Z)\subgp N_{X\ast Y}(Z)$. Thus,
	\begin{equation*}
		\Ind_{N_{X\ast Y}^{\omega}(Z)}^{N_{X\ast Y}(Z)}(\Res_{N_{X\ast Y}^{\omega}(Z)}^{\tilde{X}(\omega)\ast\tilde{Y}(\omega)}(M(X(\omega))\underset{\tilde{X}(\omega),\tilde{Y}(\omega)}{\hat{\tensor}} N(Y(\omega))))
	\end{equation*}
	is a well-defined $FN_{X\ast Y}(Z)$-module. Since $\omega$ was arbitrary, we see that the domain of $\eta_{M,N}$ is a well-defined $FN_{X\ast Y}(Z)$-module, which, moreover, is trivial source by the remarks in \ref{Brauerconsttriv} and \ref{exttensmods}(c). It is clear that $(M\hat{\tensor}_{X,Y}N)(Z)$ is a well-defined trivial source $FN_{X\ast Y}(Z)$-module; thus, it remains to check that the map $\eta_{M,N}$ is well-defined. For this, we essentially follow the same steps taken in \ref{thetaomega} above. Continue to let $\omega\in\Omega$. Consider the map
	\begin{align*}
		M(X(\omega))\times N(Y(\omega))	&\to (M\underset{X,Y}{\hat{\tensor}}N)(Z)\\
		(\Br_{X(\omega)}^M(m),\Br_{Y(\omega)}^N(n))	&\mapsto \Br_Z^{M\hat{\tensor} N}(m\tensor n),
	\end{align*}
	where $m\in M^{X(\omega)}$ and $n\in N^{Y(\omega)}$. The fact that this map is well-defined follows from an argument similar to that of \cite[3.1(h)]{Boltje_2012}. The map is clearly $F$-bilinear and balanced with respect to $k(\tilde{X}(\omega),\tilde{Y}(\omega))$, hence it descends to an $F$-linear map
	\begin{equation*}
		M(X(\omega))\underset{F[k(\tilde{X}(\omega),\tilde{Y}(\omega))]}{\tensor}N(Y(\omega))\to(M\underset{X,Y}{\hat{\tensor}}N)(Z).
	\end{equation*}
	Note that the left-hand side is the $F$-vector space underlying the extended tensor product $M(X(\omega))\hat{\tensor}_{\tilde{X}(\omega),\tilde{Y}(\omega)}N(Y(\omega))$. A straightforward calculation shows that the map above respects the action of $N_{X\ast Y}^\omega(Z)$; thus, we have a well-defined $FN_{X\ast Y}^\omega(Z)$-homomorphism
	\begin{equation*}
		\zeta_\omega:\Res_{N_{X\ast Y}^\omega(Z)}^{\tilde{X}(\omega)\ast\tilde{Y}(\omega)}(M(X(\omega))\underset{\tilde{X}(\omega),\tilde{Y}(\omega)}{\hat{\tensor}} N(Y(\omega)))\to\Res_{N_{X\ast Y}^\omega(Z)}^{N_{X\ast Y}(Z)}((M\underset{X,Y}{\hat{\tensor}}N)(Z))
	\end{equation*}
	defined by
	\begin{equation*}
		\zeta_\omega(\Br_{X(\omega)}^M(m)\tensor\Br_{Y(\omega)}^N(n))=\Br_Z^{M\hat{\tensor} N}(m\tensor n),
	\end{equation*}
	where $m\in M^{X(\omega)}$ and $n\in N^{Y(\omega)}$. By adjunction, we obtain an $FN_{X\ast Y}(Z)$-homomorphism
	\begin{equation*}
		\theta_\omega:\Ind_{N_{X\ast Y}^\omega(Z)}^{N_{X\ast Y}(Z)}(\Res_{N_{X\ast Y}^\omega(Z)}^{\tilde{X}(\omega)\ast\tilde{Y}(\omega)}(M(X(\omega))\underset{\tilde{X}(\omega),\tilde{Y}(\omega)}{\hat{\tensor}} N(Y(\omega))))\to (M\underset{X,Y}{\hat{\tensor}}N)(Z)
	\end{equation*}
	defined by
	\begin{equation*}
		\theta_\omega((g,k)\tensor(\Br_{X(\omega)}^M(m)\tensor\Br_{Y(\omega)}^N(n)))=(g,k)\cdot\Br_Z^{M\hat{\tensor} N}(m\tensor n)
	\end{equation*}
	for all $(g,k)\in N_{X\ast Y}(Z)$, $m\in M^{X(\omega)}$, and $n\in N^{Y(\omega)}$. The map $\eta_{M,N}$ in the statement of the theorem is the coproduct of the maps $\theta_\omega$, $\omega\in[T\backslash\Omega]$, and, hence, is a well-defined homomorphism of $FN_{X\ast Y}(Z)$-modules.
	
	Note that the definition of $\eta_{M,N}$ did not rely on any assumptions about the vertices of the indecomposable direct summands of $M$ and $N$. In fact, it is straightforward to check that the maps $\eta_{M,N}$ define a natural transformation $\eta$ between additive functors $\ltriv{FX}\times\ltriv{FY}\to\ltriv{FN_{X\ast Y}(Z)}$.
	
	To complete the proof we must show that $\eta_{M,N}$ is an isomorphism of $FN_{X\ast Y}(Z)$-modules if $M\in\ltriv{FX}$, $N\in\ltriv{FY}$ are such that $k(P,Q)=1$ whenever $P$ is a vertex of an indecomposable direct summand of $M$ and $Q$ is a vertex of an indecomposable direct summand of $N$. Toward this end, we may assume without loss of generality that $M$ and $N$ are indecomposable. Then by Lemma \ref{lem:forthm2}, there exist $U\in\lset{X}$, $V\in\lset{Y}$ such that $M$ is isomorphic to a direct summand of $FU$, $N$ is isomorphic to a direct summand of $FV$, and $k(X_u,Y_v)=1$ for all $u\in U$, $v\in V$. Since $\eta$ is a natural transformation between additive functors, we may assume without loss of generality that $M=FU$ and $N=FV$. But in this case, Theorem \ref{thm:1}, together with the remarks in \ref{Brauerconsttriv} and \ref{exttensmods}(c), yields that $\eta_{M,N}$ is an isomorphism. The proof is complete.
\end{proof}

\begin{remark}
	Before proceeding, we note that all of the statements in Remark \ref{rmk:natural} carry over (after the appropriate adjustments are made) into the setting of Theorem \ref{thm:2}. For example, the domain of the map $\eta_{M,N}$ in Theorem \ref{thm:2} does not depend on the choice of orbit representatives $[T\backslash\Omega]$, in the sense that any two choices lead to domains which are naturally isomorphic. Also, two different natural transformations $\eta,\eta'$ arising from two different choices of $[T\backslash\Omega]$ differ only by precomposition with this canonical isomorphism. Because these remarks are so similar to those made in Remark \ref{rmk:natural}, we leave the details to the reader.
\end{remark}

Continue to let $G$ and $H$ stand for finite groups and let $F$ be a field of positive characteristic $p$. We now assume, in addition, that $F$ is large enough for $G\times H$. Using Theorem \ref{thm:2}, we show that a $p$-permutation equivalence $\gamma$ between a block $A$ of $FG$ and a block $B$ of $FH$ induces commutative diagrams between the Grothendieck groups of trivial source modules associated to corresponding $A$- and $B$-Brauer pairs. We note that in order to apply Theorem \ref{thm:2} in this situation --- and any other previous results that involve three finite groups $G$, $H$, and $K$ --- we take $K=1$ and then identify $H$ and $G$ with $H\times K$ and $G\times K$, respectively, in the obvious way. No further comment on this convention will be made.

\begin{theorem}\label{thm:3}
	Let $A$ be a block of $FG$ with identity $e_A$, let $B$ be a block of $FH$ with identity $f_B$, and suppose that $\gamma\in T^\Delta(A,B)$ is a $p$-permutation equivalence between $A$ and $B$. Let $(\Delta(D,\phi,E),e_D\tensor f_E^\ast)$ be a maximal $\gamma$-Brauer pair. Then for any $Q\subgp E$ the diagram
	\begin{equation}\tag{$\dagger$}\label{commdiag}
	\begin{tikzcd}
		{T(B)} &&&& {T(A)} \\
		{T(FJ_Qf_Q)} &&&& {T(FI_Pe_P)}
		\arrow["{\gamma\underset{FH}{\tensor}-}", from=1-1, to=1-5]
		\arrow["{-(Q,f_Q)}"', from=1-1, to=2-1]
		\arrow["{-(P,e_P)}", from=1-5, to=2-5]
		\arrow["{\gamma(\Delta(P,\phi,Q),e_P\tensor f_Q^\ast)\underset{Y_Q,J_Q}{\hat{\tensor}}-}"', from=2-1, to=2-5]
	\end{tikzcd}
	\end{equation}
	commutes, where $P=\phi(Q)$, $e_P$ (respectively, $f_Q$) is the unique block idempotent of $FC_G(P)$ (resp., $FC_H(Q)$) such that $(P,e_P)\leq(D,e_D)$ (resp., $(Q,f_Q)\leq(E,f_E)$), $I_P=N_G(P,e_P)$, $J_Q=N_H(Q,f_Q)$, and $Y_Q=N_{G\times H}(\Delta(P,\phi,Q),e_P\tensor f_Q^\ast)$.
\end{theorem}

\begin{proof}
	Let $Q$, $P$, etc. be as in the statement. We begin by noting that each of the maps in diagram (\ref{commdiag}) are well-defined. Indeed, the vertical maps in the diagram are well-defined by \ref{gammabrauerpairs}(a), so we only need to check that the horizontal maps are also well-defined. For ease, set
	\begin{equation*}
		\gamma_Q=\gamma(\Delta(P,\phi,Q),e_P\tensor f_Q^\ast)\in T(FY_Q(e_P\tensor f_Q^\ast)).
	\end{equation*}
	Now, by \cite[Proposition 11.1]{Boltje_2020} we have $p_1(Y_Q)=I_P$, $p_2(Y_Q)=J_Q$, $k_1(Y_Q)=C_G(P)$, and $k_2(Y_Q)=C_H(Q)$. In particular, $Y_Q\ast J_Q=I_P$, so by the remarks in \ref{exttensmods}(c) the extended tensor product with $\gamma_Q$ defines a group homomorphism $\gamma_Q\hat{\tensor}_{Y_Q,J_Q}-:T_F(J_Q)\to T_F(I_P)$. This map restricts to a homomorphism $T(FJ_Qf_Q)\to T(FI_Pe_P)$ because whenever $M$ is a trivial source $FY_Q(e_P\tensor f_Q^\ast)$-module and $N$ is a trivial source $FJ_Qf_Q$-module, $e_P$ acts trivially on the $FC_G(P)$-module
	\begin{equation*}
		\Res_{C_G(P)}^{I_P}(M\underset{Y_Q,J_Q}{\hat{\tensor}}N)\iso\Res_{C_G(P)\times C_H(Q)}^{Y_Q}(M)\underset{FC_H(Q)}{\tensor}\Res_{C_H(Q)}^{J_Q}(N)
	\end{equation*}
	and, hence, the $FI_P$-module $M\hat{\tensor}_{Y_Q,J_Q}N$ belongs to the block $e_P$ of $FI_P$. We see, then, that the lower horizontal map in diagram (\ref{commdiag}) is well-defined. The well-definedness of the upper horizontal map is clear, but also follows from \ref{exttensmods}(b) and the argument just given, in the case where $Q$ is trivial.
	
	Let $\mu\in T(B)$. To complete the proof, we must show that
	\begin{equation*}
		(\gamma\underset{FH}{\tensor}\mu)(P,e_P)=\gamma_Q\underset{Y_Q,J_Q}{\hat{\tensor}}\mu(Q,f_Q).
	\end{equation*} 
	Since, by definition, the vertices of the (isomorphism classes of) indecomposable trivial source modules that support $\gamma$ are all twisted diagonal, we may apply Theorem \ref{thm:2} with $X=G\times H$, $Y=H$, and $Z=P$ to decompose $(\gamma\tensor_{FH}\mu)(P)=(\gamma\hat{\tensor}_{X,Y}\mu)(P)$ in $T_F(N_G(P))$. Before doing so, note that, in the notation of Section \ref{sec:setupthm1}, we have $\Omega=\Hom(P,H)$ and $T=N_G(P)\times H$ in this case. Now, if $[T\backslash\Omega]$ denotes a set of representatives for the $T$-orbits of $\Omega$ then by Theorem \ref{thm:2} we have
	\begin{equation*}
		(\gamma\underset{FH}{\tensor}\mu)(P)=\sum_{\omega\in[T\backslash\Omega]}\Ind_{N_G^\omega(P)}^{N_G(P)}(\Res_{N_G^\omega(P)}^{\tilde{X}(\omega)\ast\tilde{H}(\omega)}(\gamma(X(\omega))\underset{\tilde{X}(\omega),\tilde{H}(\omega)}{\hat{\tensor}}\mu(H(\omega)))).
	\end{equation*}
	(In the above and in the sequel, we continue to let $X$ stand for $G\times H$ for ease of notation.) Note that if $\omega:P\to H$ then $N_G^\omega(P)$ consists of those elements $g\in N_G(P)$ for which there exists an element $h\in H$ satisfying $\omega\circ c_g=c_h\circ\omega:P\to H$ (see \ref{N(omega)M(omega)}(a)).
	
	If $\omega:P\to H$ is such that $\gamma(X(\omega))\neq 0$ then Lemma \ref{vtxbrauer} implies that $X(\omega)$ is a twisted diagonal subgroup of $G\times H$, i.e., $\omega$ must be injective. The collection $\Inj(P,H)$ of injective homomorphisms $P\into H$ is clearly closed under the action of $T$, so we see that the sum of the previous paragraph may be reindexed over a set $[T\backslash\Inj(P,H)]$ of representatives for the orbits of $T$ on $\Inj(P,H)$. Note that if $\omega:P\into H$ is injective then we have $X(\omega)=\Delta^\omega(P)$ and $H(\omega)=\omega(P)$, the image of $P$ under $\omega$ (see \ref{twisteddiagsubg} above). Moreover, in this case we also have $N_G^\omega(P)=N_{G\times H}(\Delta^\omega(P))\ast N_H(\omega(P))$. Indeed, Lemma \ref{containment1} states that $N_G^\omega(P)$ is contained in $N_{G\times H}(\Delta^\omega(P))\ast N_H(\omega(P))$. Conversely, if $g\in N_{G\times H}(\Delta^\omega(P))\ast N_H(\omega(P))$ then there exists an element $h\in N_H(\omega(P))$ such that $(g,h)\in N_{G\times H}(\Delta^\omega(P))$. It follows easily that $g$ normalizes $p_1(\Delta^\omega(P))=P$ and $\omega\circ c_g=c_h\circ\omega:P\to H$, hence, as noted above, $g\in N_G^\omega(P)$. Noting that $N_{G\times H}(\Delta^\omega(P))=\tilde{X}(\omega)$ and $N_H(\omega(P))=\tilde{H}(\omega)$, the equality of the previous paragraph may be rewritten:
	\begin{equation}\label{eqn:1}
		(\gamma\underset{FH}{\tensor}\mu)(P)=\sum_{\omega\in[T\backslash\Inj(P,H)]}\Ind_{N_G^\omega(P)}^{N_G(P)}(\gamma(\Delta^\omega(P))\underset{\tilde{X}(\omega),\tilde{H}(\omega)}{\hat{\tensor}}\mu(\omega(P))).
	\end{equation}
	We may also assume, without loss of generality, that the set of representatives $[T\backslash\Inj(P,H)]$ indexing the sum in equation (\ref{eqn:1}) contains $\phi^{-1}:P\into H$.
	
	Let $\omega\in[T\backslash\Inj(P,H)]$. By Lemma \ref{M(P)belongsbr(e)}, $\gamma(\Delta^\omega(P))$ belongs to the $\tilde{X}(\omega)$-stable central idempotent $\br_{\Delta^\omega(P)}^{G\times H}(e_A\tensor f_B^\ast)=\br_P^G(e_A)\tensor\br_{\omega(P)}^H(f_B)^\ast$ of $F[C_G(P)\times C_H(\omega(P))]$, i.e.,
	\begin{equation*}
		\gamma(\Delta^\omega(P))\in T(F\tilde{X}(\omega)(\br_P^G(e_A)\tensor\br_{\omega(P)}^H(f_B)^\ast)).
	\end{equation*}
	So by Lemma \ref{lem:FRcor},
	\begin{equation}\label{eqn:2'}
		\gamma(\Delta^\omega(P))=\sum_{e\tensor f^\ast\in\mathcal{S}_{\omega}}\Ind_{Y_{e\tensor f^\ast}}^{\tilde{X}(\omega)}(\gamma(\Delta^\omega(P),e\tensor f^\ast))
	\end{equation}
	where $\mathcal{S}_{\omega}$ denotes a set of representatives for the $\tilde{X}(\omega)$-conjugacy classes of block idempotents $e\tensor f^\ast$ of $F[C_G(P)\times C_H(\omega(P))]$ such that $(\Delta^\omega(P),e\tensor f^\ast)$ is an $A\tensor B^\ast$-Brauer pair and where $Y_{e\tensor f^\ast}=N_{G\times H}(\Delta^\omega(P),e\tensor f^\ast)$ for each $e\tensor f^\ast\in\mathcal{S}_{\omega}$. If $\omega=\phi^{-1}$, we choose a set of representatives $\mathcal{S}_{\phi^{-1}}$ containing $e_P\tensor f_Q^\ast$.
	
	Letting $\omega\in[T\backslash\Inj(P,H)]$ be arbitrary again, we use the equality in (\ref{eqn:2'}) to rewrite the expression $\Ind_{N_G^\omega(P)}^{N_G(P)}(\gamma(\Delta^\omega(P))\hat{\tensor}_{\tilde{X}(\omega),\tilde{H}(\omega)}\mu(\omega(P)))$ which appears in equation (\ref{eqn:1}). We first note that if $e\tensor f^\ast\in\mathcal{S}_{\omega}$ then by \cite[Lemma 6.5(a)]{Boltje_2020} we have
	\begin{align*}
		&\left(\Ind_{Y_{e\tensor f^\ast}}^{\tilde{X}(\omega)}(\gamma(\Delta^\omega(P),e\tensor f^\ast))\right)\underset{\tilde{X}(\omega),\tilde{H}(\omega)}{\hat{\tensor}}\mu(\omega(P))\\
		&\qquad\qquad=\Ind_{(Y_{e\tensor f^\ast})\ast \tilde{H}(\omega)}^{N_G^\omega(P)}\left(\gamma(\Delta^\omega(P),e\tensor f^\ast)\underset{Y_{e\tensor f^\ast},\tilde{H}(\omega)}{\hat{\tensor}}\mu(\omega(P))\right).
	\end{align*}
	Therefore,
	\begin{align}\label{eqn:3'}
		&\Ind_{N_G^\omega(P)}^{N_G(P)}(\gamma(\Delta^\omega(P))\underset{\tilde{X}(\omega),\tilde{H}(\omega)}{\hat{\tensor}}\mu(\omega(P)))\\
		&\qquad\qquad=\sum_{e\tensor f^\ast\in\mathcal{S}_{\omega}}\Ind_{(Y_{e\tensor f^\ast})\ast \tilde{H}(\omega)}^{N_G(P)}\left(\gamma(\Delta^\omega(P),e\tensor f^\ast)\underset{Y_{e\tensor f^\ast},\tilde{H}(\omega)}{\hat{\tensor}}\mu(\omega(P))\right).\nonumber
	\end{align}

	Now let $\hat{e}_P$ denote the unique block idempotent of $FN_G(P)$ that covers $e_P$ (see \ref{FongReynolds} above). We claim that
	\begin{equation}\label{eqn new}
		\hat{e}_P(\gamma\underset{FH}{\tensor}\mu)(P)=\Ind_{Y_Q\ast N_H(Q)}^{N_G(P)}(\gamma_Q\underset{Y_Q,N_H(Q)}{\hat{\tensor}}\mu(Q))
	\end{equation}
	in $T(FN_G(P)\hat{e}_P)$. To see why, first note that the equalities in (\ref{eqn:1}) and (\ref{eqn:3'}) above imply that $\hat{e}_P(\gamma\underset{FH}{\tensor}\mu)(P)$ is equal to
	\begin{equation}\label{eqn:4'}
		\sum_{\substack{\omega\in[T\backslash\Inj(P,H)]\\e\tensor f^\ast\in\mathcal{S}_{\omega}}}\hat{e}_P\Ind_{(Y_{e\tensor f^\ast})\ast \tilde{H}(\omega)}^{N_G(P)}\left(\gamma(\Delta^\omega(P),e\tensor f^\ast)\underset{Y_{e\tensor f^\ast},\tilde{H}(\omega)}{\hat{\tensor}}\mu(\omega(P))\right).
	\end{equation}
	Suppose that $\omega\in[T\backslash\Inj(P,H)]$ and $e\tensor f^\ast\in\mathcal{S}_{\omega}$ are such that the $(\omega,e\tensor f^\ast)$-term of the sum in (\ref{eqn:4'}) is nonzero. Then we must have $\gamma(\Delta^\omega(P),e\tensor f^\ast)\neq 0$, i.e., $(\Delta^\omega(P),e\tensor f^\ast)$ must be a $\gamma$-Brauer pair. Furthermore, the block idempotent $e$ of $FC_G(P)$ must be $N_G(P)$-conjugate to $e_P$. Indeed, it is easy to see that
	\begin{equation*}
		PC_G(P)\subgp(Y_{e\tensor f^\ast})\ast\tilde{H}(\omega)\subgp N_G(P,e),
	\end{equation*}
	so that $e$ is a block idempotent of $F[(Y_{e\tensor f^\ast})\ast\tilde{H}(\omega)]$ by Lemma \ref{lem:forthm3}. It is also clear that
	\begin{equation*}
		\gamma(\Delta^\omega(P),e\tensor f^\ast)\underset{Y_{e\tensor f^\ast},\tilde{H}(\omega)}{\hat{\tensor}}\mu(\omega(P))\in T(F[(Y_{e\tensor f^\ast})\ast\tilde{H}(\omega)]e).
	\end{equation*}
	Therefore, if $\hat{e}$ denotes the unique block idempotent of $FN_G(P)$ covering $e$ then 
	\begin{equation*}
		\Ind_{(Y_{e\tensor f^\ast})\ast \tilde{H}(\omega)}^{N_G(P)}\left(\gamma(\Delta^\omega(P),e\tensor f^\ast)\underset{Y_{e\tensor f^\ast},\tilde{H}(\omega)}{\hat{\tensor}}\mu(\omega(P))\right)
	\end{equation*}
	is an element of $T(FN_G(P)\hat{e})$, again by Lemma \ref{lem:forthm3}. If the $\hat{e}_P$-component of this expression is nonzero then we must have $\hat{e}=\hat{e}_P$, which can only occur if $e$ is $N_G(P)$-conjugate to $e_P$, as needed.
	
	Let $g\in N_G(P)$ be such that ${}^ge_P=e$. Since the set of $\gamma$-Brauer pairs is stable under conjugation by elements of $G\times H$ (see the remarks in \ref{gammabrauerpairs}(b)),
	\begin{equation*}
		{}^{(g^{-1},1)}(\Delta^\omega(P),e\tensor f^\ast)=(\Delta^{\omega c_g}(P),e_P\tensor f^\ast)
	\end{equation*}
	is a $\gamma$-Brauer pair. Since all maximal $\gamma$-Brauer pairs are $G\times H$-conjugate, there exists a pair $(g_1,h_1)\in G\times H$ such that 
	\begin{equation*}
		{}^{(g_1,h_1)}(\Delta^{\omega c_g}(P),e_P\tensor f^\ast)\leq (\Delta(D,\phi,E),e_D\tensor f_E^\ast).
	\end{equation*}
	This containment implies that ${}^{g_1}(P,e_P)\leq(D,e_D)$ in $\BP(A)$\footnote{In other words, conjugation with $g_1$ defines an isomorphism between $P$ and ${}^{g_1}P$ in the fusion system on $D$ associated to the block $A$.} and ${}^{h_1}(\omega(P),f)\leq(E,f_E)$ in $\BP(B)$. By \cite[Theorem 11.2]{Boltje_2020} there exists an element $h\in H$ such that ${}^h(Q,f_Q)={}^{h_1}(\omega(P),f)$ and for which the diagram below commutes:
	\[\begin{tikzcd}
		P & {{}^{g_1}P} \\
		Q & {{}^{h_1}\omega(P)}
		\arrow["{c_{g_1}}", from=1-1, to=1-2]
		\arrow["\phi", from=2-1, to=1-1]
		\arrow["{c_h}"', from=2-1, to=2-2]
		\arrow["\phi"', from=2-2, to=1-2]
	\end{tikzcd}\]
	Now, since ${}^{(g_1,h_1)}\Delta^{\omega c_g}(P)\subgp \Delta(D,\phi,E)$, the restriction of the isomorphism $\phi$ to ${}^{h_1}\omega(P)=p_2({}^{(g_1,h_1)}\Delta^{\omega c_g}(P))$ is equal to $c_{g_1g^{-1}}\omega^{-1}c_{h_1}^{-1}:{}^{h_1}\omega(P)\to {}^{g_1}P$. It follows that 
	\begin{equation*}
		\phi^{-1}\circ c_{g^{-1}}=c_{h^{-1}h_1}\circ\omega:P\to Q.
	\end{equation*}
	Since $g\in N_G(P)$ and $T=N_G(P)\times H$ the above shows that $\omega$ is $T$-conjugate to $\phi^{-1}:P\into H$. Since the set of representatives $[T\backslash\Inj(P,H)]$ contains $\phi^{-1}$ it follows that $\omega=\phi^{-1}$. But then the equalities above yield that $(g,h_1^{-1}h)\in N_{G\times H}(\Delta(P,\phi,Q))$ and ${}^{(g,h_1^{-1}h)}(e_P\tensor f_Q^\ast)=e\tensor f^\ast$. Thus, $e\tensor f^\ast\in\mathcal{S}_{\phi^{-1}}$ is conjugate by an element of $\tilde{X}(\phi^{-1})=N_{G\times H}(\Delta(P,\phi,Q))$ to $e_P\tensor f_Q^\ast$. By our choice of $\mathcal{S}_{\phi^{-1}}$ we must have $e\tensor f^\ast=e_P\tensor f_Q^\ast$. 
	
	We have shown that the only possible nonzero term in the sum of (\ref{eqn:4'}) above is the term indexed by $\phi^{-1}\in [T\backslash\Inj(P,H)]$ and $e_P\tensor f_Q^\ast\in \mathcal{S}_{\phi^{-1}}$. Since the sum in (\ref{eqn:4'}) is equal to $\hat{e}_P(\gamma\underset{FH}{\tensor}\mu)(P)$, and since $Y_{e_P\tensor f_Q^\ast}=Y_Q$ and $\tilde{H}(\phi^{-1})=N_H(Q)$, we find that
	\begin{align*}
		\hat{e}_P(\gamma\underset{FH}{\tensor}\mu)(P)	&=\hat{e}_P\Ind_{Y_Q\ast N_H(Q)}^{N_G(P)}(\gamma_Q\underset{Y_Q,N_H(Q)}{\hat{\tensor}}\mu(Q))\\
		&=\Ind_{Y_Q\ast N_H(Q)}^{N_G(P)}(\gamma_Q\underset{Y_Q,N_H(Q)}{\hat{\tensor}}\mu(Q)).
	\end{align*}
	This completes the proof of the claim, namely Equation~(\ref{eqn new}).
	
	Now, it is easy to verify that $Y_Q\ast N_H(Q)=I_P$, so the equality above can be written:
	\begin{equation*}
		\hat{e}_P(\gamma\underset{FH}{\tensor}\mu)(P)=\Ind_{I_P}^{N_G(P)}(\gamma_Q\underset{Y_Q,N_H(Q)}{\hat{\tensor}}\mu(Q)).
	\end{equation*}
	Since $\gamma_Q\in T(FY_Q(e_P\tensor f_Q^\ast))$, if $\hat{f}_Q$ denotes the unique block idempotent of $FN_H(Q)$ that covers $f_Q$ then we have
	\begin{equation*}
		\gamma_Q\underset{Y_Q,N_H(Q)}{\hat{\tensor}}\mu(Q)=\gamma_Q\underset{Y_Q,N_H(Q)}{\hat{\tensor}}\hat{f}_Q\mu(Q)=\gamma_Q\underset{Y_Q,N_H(Q)}{\hat{\tensor}}\Ind_{J_Q}^{N_H(Q)}(\mu(Q,f_Q)),
	\end{equation*}
	the final equality following from Lemma \ref{lem:FRcor}. But it is clear from the definition of the extended tensor product that
	\begin{equation*}
		\gamma_Q\underset{Y_Q,N_H(Q)}{\hat{\tensor}}\Ind_{J_Q}^{N_H(Q)}(\mu(Q,f_Q))=\gamma_Q\underset{Y_Q,J_Q}{\hat{\tensor}}\mu(Q,f_Q).
	\end{equation*}
	It follows that 
	\begin{equation*}
		\hat{e}_P(\gamma\underset{FH}{\tensor}\mu)(P)=\Ind_{I_P}^{N_G(P)}(\gamma_Q\underset{Y_Q,J_Q}{\hat{\tensor}}\mu(Q,f_Q))
	\end{equation*}
	in $T(FN_G(P)\hat{e}_P)$. Applying the Fong-Reynolds correspondence (\ref{FongReynolds}) allows us to conclude that
	\begin{equation*}
		(\gamma\underset{FH}{\tensor}\mu)(P,e_P)=\gamma_Q\underset{Y_Q,J_Q}{\hat{\tensor}}\mu(Q,f_Q)
	\end{equation*}
	in $T(FI_Pe_P)$. The proof is complete.
\end{proof}

\begin{remark}\label{rmk:isotypydiagrams}
	The commutative diagrams of Theorem \ref{thm:3} can be considered as liftings of the commutative diagrams that appear in Brou\'{e}'s concept of an \textit{isotypy} between blocks (cf. \cite{Broue_1990}) to the level of trivial source modules. Let $(\K,\OO,F)$ be a $p$-modular system large enough for $G\times H$, let $A=\OO Ge_A$ be a block of $\OO G$ and let $B=\OO Hf_B$ be a block of $\OO H$. Many of the concepts recalled in Section \ref{sec:trivsourcemodsandbrauerpairs}, which were there defined over $F$, may just as well be defined over $\OO$. For example, it is known that each trivial source $F\tensor_\OO A$-module has a unique (up to isomorphism) lift to a trivial source $A$-module (see, e.g., \cite[Theorem 5.10.2(iv)]{Linckelmann_2018}). In particular, $T(A)\iso T(F\tensor_\OO A)$ and the commutative diagrams of Theorem \ref{thm:3} also hold over $\OO$. Let $\gamma\in T^\Delta(A,B)$ be a $p$-permutation equivalence and let $(\Delta(D,\phi,E),e_D\tensor f_E^\ast)$ be a maximal $\gamma$-Brauer pair. Let $v\in E$ and let $f_{\gp{v}}$ denote the unique block idempotent of $\OO C_H(v)$ such that $(\gp{v},f_{\gp{v}})\leq(E,f_E)$ is a containment of $B$-Brauer pairs. Set $u=\phi(v)$ and let $e_{\gp{u}}$ denote the unique block idempotent of $\OO C_G(u)$ such that $(\gp{u},e_{\gp{u}})\leq(D,e_D)$ is a containment of $A$-Brauer pairs. Let $(\mu_Q)_{Q\subgp E}$ be the isotypy obtained from $\gamma$ as in \cite[Theorem 15.4]{Boltje_2020}. Consider the diagram
	\[\begin{tikzcd}
		{T(B)} &&& {T(A)} \\
		& {CF(H,B)} & {CF(G,A)} \\
		& {CF(C_H(v),f_{\gp{v}})} & {CF(C_G(u),e_{\gp{u}})} \\
		{T(\OO J_{\gp{v}}f_{\gp{v}})} &&& {T(\OO I_{\gp{u}}e_{\gp{u}})}
		\arrow["{\gamma\underset{\OO H}{\tensor}-}", from=1-1, to=1-4]
		\arrow[from=1-1, to=2-2]
		\arrow["{-(\gp{v},f_{\gp{v}})}", from=1-1, to=4-1]
		\arrow[from=1-4, to=2-3]
		\arrow["{-(\gp{u},e_{\gp{u}})}"', from=1-4, to=4-4]
		\arrow["{\mu_{\set{1}}\underset{\K H}{\tensor}-}", from=2-2, to=2-3]
		\arrow["{d_H^{v,f_{\gp{v}}}}"', from=2-2, to=3-2]
		\arrow["{d_G^{u,e_{\gp{u}}}}", from=2-3, to=3-3]
		\arrow["{\mu_{\gp{v}}\underset{\K C_H(v)}{\tensor}-}", from=3-2, to=3-3]
		\arrow[from=4-1, to=3-2]
		\arrow["{\gamma(\Delta(\gp{u},\phi,\gp{v}),e_{\gp{u}}\tensor f_{\gp{v}}^\ast)\underset{Y_{\gp{v}},J_{\gp{v}}}{\hat{\tensor}}-}"', from=4-1, to=4-4]
		\arrow[from=4-4, to=3-3]
	\end{tikzcd}\]
	Within the inner square of this diagram, the vertices are spaces of $\K$-valued class functions that belong to the blocks specified, the vertical arrows are generalized decomposition maps (cf. \cite[Definition 6.13.3]{Linckelmann_2018_2}), and the horizontal maps are perfect isometries arising from the isotypy $(\mu_Q)_{Q\subgp E}$ (cf. \cite{Broue_1990}, \cite[Chapter 9.2]{Linckelmann_2018_2}). The upper diagonal maps are induced by ``taking characters,'' while the lower diagonal maps are induced by taking characters and then restricting to centralizers.
	
	The diagram above commutes. Indeed, the inner square commutes because $(\mu_Q)_{Q\subgp E}$ is an isotypy, and the outer square commutes by Theorem \ref{thm:3}. Now, by definition (see \cite[Theorem 15.4]{Boltje_2020}), $\mu_{\set{1}}$ is the character of $\gamma$ and $\mu_{\gp{v}}$ is the character of the restriction of an $\OO$-lift of $\gamma(\Delta(\gp{u},\phi,\gp{v}),e_{\gp{u}}\tensor f_{\gp{v}}^\ast)$ to $C_G(u)\times C_H(v)$. Thus, the upper face of the diagram clearly commutes, while the lower face commutes thanks to the remarks in \ref{exttensmods}(b). Finally, the fact that the left and right faces commute follows from \cite[Corollary 2.6]{Boltje_2008}. 
	
	The diagram above establishes the connection between the commutative diagrams appearing in Theorem \ref{thm:3} and those of an isotypy between blocks. While the diagrams of Theorem \ref{thm:3} exist at the level of trivial source modules rather than at the level of class functions, one has such a diagram for each subgroup $Q$ of the defect group $E$ of the block $B$. This is in contrast to the situation of an isotypy, in which case one has a commutative diagram for each \textit{cyclic} subgroup of $E$.
\end{remark}


\bibliographystyle{plain}
\bibliography{../../../../Bibliography/bibliography}

\begin{thebibliography}{10}

\bibitem{Alperin_1979}
J.~L. Alperin and M.~Brou\'{e}.
\newblock Local methods in block theory.
\newblock {\em Ann. of Math.}, 110(1):143--157, 1979.

\bibitem{Aschbacher_2011}
M.~Aschbacher, R.~Kessar, and B.~Oliver.
\newblock {\em Fusion Systems in Algebra and Topology}.
\newblock London Mathematical Society Lecture Note Series. Cambridge University
  Press, 2011.

\bibitem{Benson_1991_vol1}
D.~J. Benson.
\newblock {\em Representations and Cohomology I: Basic Representation Theory of
  Finite Groups and Associative Algebras}, volume~30 of {\em Cambridge Studies
  in Advanced Mathematics}.
\newblock Cambridge University Press, 1991.

\bibitem{Boltje_2012}
R.~Boltje and S.~Danz.
\newblock A ghost ring for the left-free double burnside ring and an
  application to fusion systems.
\newblock {\em Advances in Mathematics}, 229:1688--1733, 2012.

\bibitem{Boltje_2020}
R.~Boltje and P.~Perepelitsky.
\newblock $p$-{P}ermutation equivalences between blocks of group algebras.
\newblock {\em Journal of Algebra}, 664:815--887, 2025.

\bibitem{Boltje_2008}
R.~Boltje and B.~Xu.
\newblock On $p$-permutation equivalences: between {R}ickard equivalences and
  isotypies.
\newblock {\em Transactions of the American Mathematical Society},
  360(10):5067--5087, 2008.

\bibitem{Bouc_2010}
S.~Bouc.
\newblock {\em Biset Functors for Finite Groups}, volume 1990 of {\em Lecture
  Notes in Mathematics}.
\newblock Springer-Verlag Berlin Heidelberg, 2010.

\bibitem{Bouc_2010.1}
S.~Bouc.
\newblock Bisets as categories and tensor product of induced bimodules.
\newblock {\em Applied Categorical Structures}, 18(5):517--521, 2010.
\newblock doi: 10.1007/s10485-008-9180-1.

\bibitem{Broue_1985}
M.~Broué.
\newblock On {S}cott modules and {$p$}-permutation modules: an approach through
  the {B}rauer morphism.
\newblock {\em Proceedings of the American Mathematical Society},
  93(3):401--408, 1985.

\bibitem{Broue_1990}
M.~Brou\'{e}.
\newblock Isom\'etries parfaites, types de blocs, cat\'egories d\'eriv\'ees.
\newblock In {\em Repr\'esentations lin\'eaires des groupes finis - Luminy,
  16-21 mai 1988}, number 181-182 in Ast\'erisque, pages 61--92. Soci\'et\'e
  math\'ematique de France, 1990.

\bibitem{Linckelmann_2018_2}
M.~Linckelmann.
\newblock {\em The Block Theory of Finite Group Algebras}, volume~2 of {\em
  London Mathematical Society Student Texts}.
\newblock Cambridge University Press, 2018.

\bibitem{Linckelmann_2018}
M.~Linckelmann.
\newblock {\em The Block Theory of Finite Group Algebras}, volume~1 of {\em
  London Mathematical Society Student Texts}.
\newblock Cambridge University Press, 2018.

\bibitem{Nagao_1989}
H.~Nagao and Y.~Tsushima.
\newblock {\em Representations of Finite Groups}.
\newblock Academic Press, San Diego, 1989.

\end{thebibliography}

\end{document}